\documentclass[letterpaper, 11pt]{article}
\usepackage[margin=3.4cm]{geometry}
\usepackage[round]{natbib}

\RequirePackage[colorlinks,citecolor=blue,urlcolor=blue, linkcolor=blue]{hyperref}
\usepackage{comment}
\usepackage{lineno}
\usepackage{amsmath, amsthm, amssymb, mathtools, bm, bbm, mathrsfs}
\usepackage{graphicx} 
\usepackage{tikz}
\usetikzlibrary{arrows, decorations.pathreplacing, positioning, shapes,arrows.meta}
\usepackage[capitalise]{cleveref}

\Crefname{assumption}{Assumption}{Assumptions}
\usepackage{chngcntr}
\usepackage{apptools}
\usepackage{caption}
\usepackage{subcaption}
\usepackage{booktabs, multirow}
\usepackage{enumerate}
\usepackage{algpseudocode, subcaption}
\usepackage{enumitem}

\makeatletter
\def\cctext#1{\expandafter\@cctext\csname c@#1\endcsname}
\def\@cctext#1{\ifcase#1\or \textbf{Method 1}\or \textbf{Method 2}\or \textbf{Method 3}\or \textbf{Method 4}\fi}
\makeatother
\AddEnumerateCounter{\cctext}{\@cctext}{Method}

\usepackage[ruled,linesnumbered]{algorithm2e}
\SetKwInOut{Input}{input}
\SetKwInOut{Output}{output}
\SetKw{Throw}{throw}

\newcommand*\dd{\mathop{}\!\mathrm{d}}
\DeclareMathOperator{\E}{\mathbb{E}}

\DeclareMathOperator{\DAG}{DAG}

\DeclareMathOperator{\MPDAG}{MPDAG}

\DeclareMathOperator{\Cov}{cov}
\DeclareMathOperator{\Var}{var}

\DeclareMathOperator{\Do}{do}

\DeclareMathOperator{\PossDe}{PossDe}

\DeclareMathOperator{\An}{An}

\DeclareMathOperator{\Pa}{Pa}
\DeclareMathOperator{\pa}{Pa}

\newcommand{\mb}[1]{\mathbf{#1}}

\newcommand{\mpdag}{MPDAG}

\newcommand{\vars}[1][V]{\mathbf{#1}}

\newcommand{\g}[1][G]{\mathcal{#1}}

\newcommand{\pstar}[1][p]{{#1}^{*}}

\definecolor{blue-violet}{rgb}{0.54, 0.17, 0.89}
\definecolor{antiquefuchsia}{rgb}{0.57, 0.36, 0.51}
\definecolor{amethyst}{rgb}{0.6, 0.4, 0.8}
\definecolor{blue-violet}{rgb}{0.54, 0.17, 0.89}
\definecolor{ao}{rgb}{0.0, 0.5, 0.0}
\definecolor{blue(ncs)}{rgb}{0.0, 0.53, 0.74}
\definecolor{dgreen}{rgb}{0.12, 0.3, 0.17}
\definecolor{cadmiumgreen}{rgb}{0.0, 0.42, 0.24}
\definecolor{darkolivegreen}{rgb}{0.33, 0.42, 0.18}
\definecolor{dartmouthgreen}{rgb}{0.05, 0.5, 0.06}

\theoremstyle{plain}
\newtheorem{theorem}{Theorem}
\newtheorem{corollary}[theorem]{Corollary}
\newtheorem{lemma}[theorem]{Lemma}

\theoremstyle{definition}
\newtheorem{example}[theorem]{Example}

 \newenvironment{proofof}[1][]{\begin{trivlist}
 \item[\hskip \labelsep {\bfseries Proof of #1.}]}{\hfill{}$\square$\end{trivlist}}

\theoremstyle{remark}

\let \epsilon \varepsilon

\title{Minimal Enumeration of All Possible Total Effects in a Markov Equivalence Class}

\author{ {\bf F.\ Richard Guo} \\Department of Statistics \\ University of Washington\\
\tt{ricguo@uw.edu}\\ \and 
{\bf Emilija Perkovi\'c} \\Department of Statistics \\ University of Washington\\ \tt{perkovic@uw.edu} \\
}

\makeatletter\begin{document}
\maketitle

\begin{abstract}
In observational studies, when a total causal effect of interest is not identified, the set of all possible effects can be reported instead. This typically occurs when the underlying causal DAG is only known up to a Markov equivalence class, or a refinement thereof due to background knowledge. As such, the class of possible causal DAGs is represented by a maximally oriented partially directed acyclic graph (MPDAG), which contains both directed and undirected edges. We characterize the minimal additional edge orientations required to identify a given total effect. A recursive algorithm is then developed to enumerate subclasses of DAGs, such that the total effect in each subclass is identified as a distinct functional of the observed distribution. This resolves an issue with existing methods, which often report possible total effects with duplicates, namely those that are numerically distinct due to sampling variability but are in fact causally identical.
\end{abstract}

\section{Introduction}

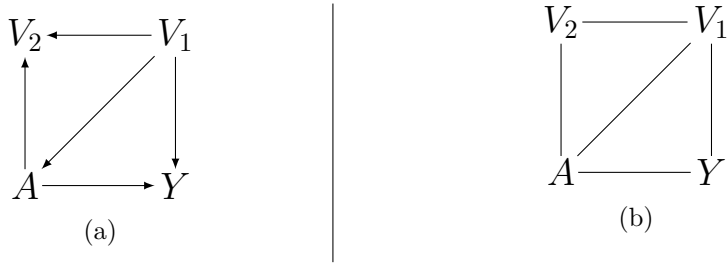
\begin{figure}[!t]
\begin{subfigure}{.4\columnwidth}
  \centering
\begin{tikzpicture}[->,>=latex,shorten >=1pt,auto,node distance=0.8cm,scale=1,transform shape]
  \tikzstyle{state}=[inner sep=1pt, minimum size=12pt]

\node[state] (Xia) at (0,0) {\Large $A$};
  \node[state] (Xka) at (0,2) {\Large $V_2$};
  \node[state] (Xja) at (2,0) {\Large  $Y$};
  \node[state] (Xsa) at (2,2) {\Large  $V_1$};

  \draw (Xka) edge [<-] (Xsa);
\draw 		(Xia) edge  [->] (Xka);
 \draw   	(Xja) edge [<-] (Xsa);
 \draw    	(Xia) edge [->] (Xja);
\draw    	(Xia) edge [<-] (Xsa);
\end{tikzpicture}
\caption{}
  \label{exdag-dag}
\end{subfigure}
\vrule
\hspace{0.2cm}
\begin{subfigure}{.5\columnwidth}
\centering
\begin{tikzpicture}[->,>=latex,shorten >=1pt,auto,node distance=0.8cm,scale=1,transform shape]
  \tikzstyle{state}=[inner sep=1pt, minimum size=12pt]

\node[state] (Xia) at (0,0) {\Large $A$};
  \node[state] (Xka) at (0,2) {\Large $V_2$};
  \node[state] (Xja) at (2,0) {\Large  $Y$};
  \node[state] (Xsa) at (2,2) {\Large  $V_1$};

  \draw (Xka) edge [-] (Xsa);
\draw 		(Xia) edge  [-] (Xka);
 \draw   	(Xja) edge [-] (Xsa);
 \draw    	(Xia) edge [-] (Xja);
\draw    	(Xia) edge [-] (Xsa);
\end{tikzpicture}
\caption{}
\label{exdag-cpdag}
\hspace{.5cm}
\end{subfigure}
\caption{(a) DAG  $\g[D]$, (b) CPDAG $\g[C]$.} 
\label{fig:example_DAG}
\end{figure}

We consider identifying total causal effects (``total effects'' or simply ``effects'' throughout) from causal graphs that can be learned from observational data and background knowledge, under the assumption of no latent variables. The full knowledge of the causal system is typically represented by a directed acyclic graph (DAG) \citep{pearl2009causality}. 
\cref{fig:example_DAG}(a) shows an example DAG $\g[D]$. Each node $u \in \{A,Y,V_1,V_2\}$ in $\g[D]$ represents a random variable $X_{u}$ in a random vector $X = (X_A, X_Y,X_{V_1}, X_{V_2})$. Each edge in $\g[D]$ represents a direct causal relationship between two variables.

Given a causal DAG, every total effect can be identified and hence consistently estimated from observational data \citep{robins1986new,pearl1995causal,pearl1995probabilistic,galles1995testing}.
In general, however, one cannot learn a causal DAG from observational data. Instead, under the assumption of no latent variables, one can learn a Markov equivalence class of DAGs that can give rise to the observed distribution. A Markov equivalence class is uniquely represented by a completed partially directed acyclic graph (CPDAG), also known as an essential graph \citep{meek1995causal,andersson1997characterization,spirtes2000causation,Chickering02}. Within the equivalence class, one DAG should not be preferred over another based on observational data. \cref{fig:example_DAG}(b) shows the CPDAG $\g[C]$ that represents $\g[D]$.

\begin{figure}[!t]
\tikzstyle{every edge}=[draw,>=stealth',->]
\newcommand\dagvariant[1]{\begin{tikzpicture}[xscale=.5,yscale=0.5]
\node (a) at (0,0) {};
\node (d) at (0,2) {};
\node (b) at (2,0) {};
\node (c) at (2,2) {};
\draw (a) edge [-] (b);
\draw (b) edge [-] (c);
\draw (d) edge [-] (a);
\draw (d) edge [-] (c);
\draw (a) edge [-] (c);
\draw #1;
\end{tikzpicture}}

\tikzstyle{every edge}=[draw,>=stealth',->]
\newcommand\dagempty[1]{\begin{tikzpicture}[xscale=.5,yscale=0.5]
\node (a) at (0,0) {};
\node (d) at (0,.5) {};
\node (b) at (2,0) {};
\node (c) at (2,.5) {};
\draw #1;
\end{tikzpicture}}

\centering
\begin{subfigure}{.4\columnwidth}
  \centering
\begin{tikzpicture}[->,>=latex,shorten >=1pt,auto,node distance=0.8cm,scale=1,transform shape]
  \tikzstyle{state}=[inner sep=1pt, minimum size=12pt]

\node[state] (Xia) at (0,0) {\Large $A$};
  \node[state] (Xka) at (0,2) {\Large $V_2$};
  \node[state] (Xja) at (2,0) {\Large  $Y$};
  \node[state] (Xsa) at (2,2) {\Large  $V_1$};

  \draw (Xka) edge [-] (Xsa);
\draw 		(Xia) edge  [-] (Xka);
 \draw   	(Xja) edge [-] (Xsa);
 \draw    	(Xia) edge [->,line width=1pt] (Xja);
\draw    	(Xia) edge [-] (Xsa);
\end{tikzpicture}
\caption{}
  \label{mpdag}
\end{subfigure}
\vrule
\hspace{0.2cm}
\begin{subfigure}{.5\columnwidth}
\begin{center}
\dagvariant{
(a)  edge [->,line width=1pt]  (b)
(d)  edge [->,line width=1pt]  (a)
(a) edge [->,line width=1pt] (c)
(d) edge [->] (c)
(c) edge [->] (b)
}\\
\end{center}
\begin{center}
\dagvariant{
(a)  edge [->,line width=1pt]  (b)
(a)  edge [->,line width=1pt]  (d)
(a) edge [->,line width=1pt] (c)
(d) edge [->] (c)
(c) edge [->] (b)
}
\dagvariant{
(a)  edge [->,line width=1pt]  (b)
(a)  edge [->,line width=1pt]  (d)
(a) edge [->,line width=1pt] (c)
(c) edge [->] (d)
(c) edge [->] (b)
}
\dagvariant{
(a)  edge [->,line width=1pt]  (b)
(a)  edge [->,line width=1pt]  (d)
(a) edge [->,line width=1pt] (c)
(c) edge [->] (d)
(b) edge [->] (c)
}
\end{center}
\begin{center}
\dagvariant{
(a)  edge [->,line width=1pt]  (b)
(a)  edge [->,line width=1pt,color=blue]  (d)
(c) edge [->,line width=1pt,color=blue] (a)
(c) edge [->] (d)
(c) edge [->] (b)
}
\dagvariant{
(a)  edge [->,line width=1pt]  (b)
(d)  edge [->,line width=1pt,color=purple]  (a)
(c) edge [->,line width=1pt,color=purple] (a)
(d) edge [->] (c)
(c) edge [->] (b)
}
\dagvariant{
(a)  edge [->,line width=1pt]  (b)
(d)  edge [->,line width=1pt,color=purple]  (a)
(c) edge [->,line width=1pt,color=purple] (a)
(c) edge [->] (d)
(c) edge [->] (b)
}
\end{center}
\caption{}
\label{alldags-mpdag}
\hspace{.5cm}
\end{subfigure}
\begin{subfigure}{.4\columnwidth}
\begin{center}
\dagvariant{
(a)  edge [->,line width=1pt]  (b)
(d)  edge [->,line width=1pt]  (a)
(a) edge [->,line width=1pt] (c)
(d) edge [->] (c)
(c) edge [->] (b)
}
\dagvariant{
(a)  edge [->,line width=1pt]  (b)
(a)  edge [->,line width=1pt]  (d)
(a) edge [->,line width=1pt] (c)
}
\end{center}
\begin{center}
\dagvariant{
(a)  edge [->,line width=1pt]  (b)
(a)  edge [->,line width=1pt,color=blue]  (d)
(c) edge [->,line width=1pt,color=blue] (a)
(c) edge [->] (d)
(c) edge [->] (b)
}
\dagvariant{
(a)  edge [->,line width=1pt]  (b)
(d)  edge [->,line width=1pt,color=purple]  (a)
(c) edge [->,line width=1pt,color=purple] (a)
(c) edge [->] (b)
}
\end{center}
\caption{}
\label{alldags-mpdag2}
\end{subfigure}
\vrule
\hspace{0.2cm}
\begin{subfigure}{.5\columnwidth}
\begin{center}
\dagvariant{
(a)  edge [->,line width=1pt]  (b)
(d)  edge [->,line width=1pt]  (a)
(a) edge [->,line width=1pt] (c)
(d) edge [->] (c)
(c) edge [->] (b)
}
\dagvariant{
(a)  edge [->,line width=1pt]  (b)
(a)  edge [->,line width=1pt]  (d)
(a) edge [->,line width=1pt] (c)
}
\end{center}
\begin{center}
\dagvariant{
(a)  edge [->,line width=1.3pt]  (b)
(c) edge [->,line width=1.3pt,color=blue-violet] (a)
(c) edge [->] (b)
}
\end{center}
\caption{}
\label{alldags-mpdag3}
\end{subfigure}
\caption{\cref{exam:point}. (a) MPDAG $\g$, (b) all DAGs represented by $\g$, (c) all MPDAGs with distinct parent sets of $A$ represented by $\g$, (d) all MPDAGs with distinct causal identification formulas for $f(x_Y|\Do(x_A))$ represented by $\g$.} 
\label{fig:example_IDA}
\end{figure}
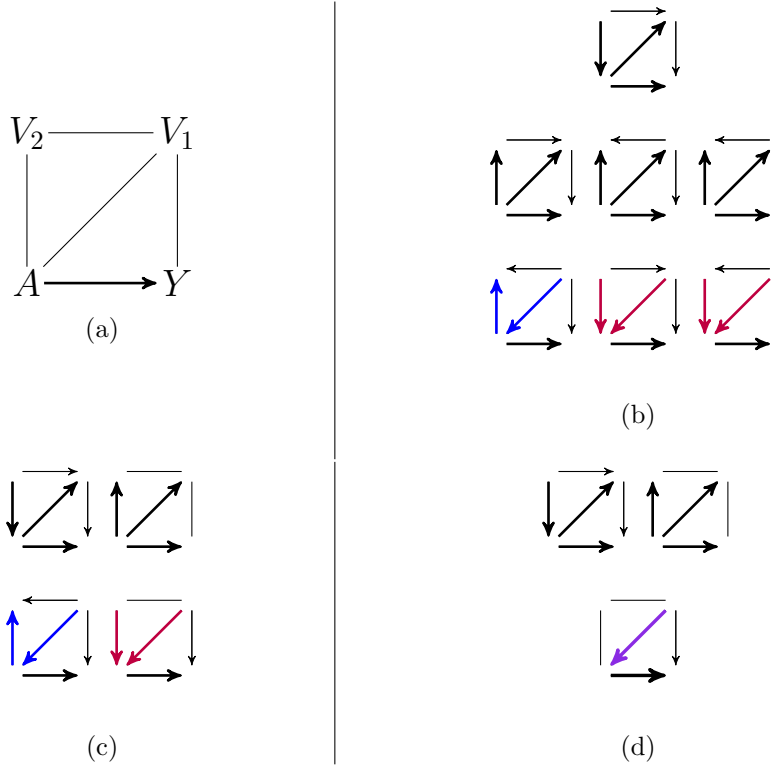

Often, we may have additional background knowledge on the underlying causal system. For example, we may know that $A$ temporally precedes $Y$ and therefore determine (or reveal) the edge orientation $A \rightarrow Y$ in CPDAG $\g[C]$.  This results in a maximally oriented partially directed acyclic graph (MPDAG) $\g$, drawn in \cref{fig:example_IDA}(a). MPDAGs are a class of graphs that subsumes both CPDAGs and DAGs. They are obtained by (optionally) adding edge orientations to a CPDAG and completing the orientation rules of \cite{meek1995causal} (see \cref{fig:orientationRules}). As such, the class of DAGs represented by an MPDAG is a refinement of the corresponding Markov equivalence class. For example, the class of DAGs represented by $\g$ is drawn in \cref{fig:example_IDA}(b), which consists of all DAGs in the Markov equivalence class represented by $\g[C]$ that have $A \rightarrow Y$.

The background knowledge of pairwise causal relationships of this type can be derived from field expertise. Moreover, other types of background knowledge, such as tiered orderings \citep{tetrad1998}, non-ancestral background knowledge \citep{fang2020ida}, knowledge derived from  experimental data \citep{hauserBuehlmann12,wang2017permutation} as well as certain model restrictions \citep{hoyer08,ernestroth2016} can also be used to obtain an MPDAG.

A total effect is identified given an equivalence class of DAGs if it can be expressed as a functional of the observed distribution, which is the same for all DAGs in the equivalence class; see \cref{sec:prelim} for the definition. 
Recently, \cite{perkovic20} gave a necessary and sufficient graphical condition for identifying a total effect given an MPDAG  (Theorem \ref{thm:id-criterion}).
When the condition fails, the effect of interest cannot be identified, such as the effect of $A$ on $Y$ given $\g$ in \cref{fig:example_IDA}(a).
In such cases, the observational data can still be informative if one identifies a finite set that contains the true effect. 
To do so, one can enumerate all DAGs represented by $\g$ and estimate the total effect under each, obtaining a set of estimated \emph{possible} total effects. For instance, for the class of DAGs in \cref{fig:example_IDA}(b), 7 possible total effects can be reported. 

However, there are two drawbacks to this approach. First, enumerating all DAGs in a Markov equivalence class (or a refinement thereof) is computationally prohibitive unless one has only a few variables. For example, the complete CPDAG of $p$ variables contains $p!$ DAGs; see also \citet{gillispie2002size,steinsky2013enumeration}. Second, the number of distinct possible effects can be much smaller than the size of the equivalence class. For example, the effect of $A$ on $Y$ is the same for the three DAGs listed in the second row of \cref{fig:example_IDA}(b). That being said, depending on the estimator applied to each DAG, one may obtain three estimates that only \emph{look} different in finite samples, which are in fact different estimators for the same possible effect. These statistical \emph{duplicates} are undesirable as it undermines the interpretability of the estimated set.
Hence, to save computation time and to deliver causally informative estimates, one should instead enumerate all possible effects that are distinct, or in other words, \emph{minimally}.

Recent works on this topic include the ``intervention calculus when the DAG is absent'' (IDA) algorithms and joint-IDA algorithms of \cite{maathuis2009estimating}, \cite{nandy2017estimating}, \cite{perkovic17}, \cite{witte2020efficient} and \citet{fang2020ida}. Given an MPDAG $\g$, these methods enumerate a set of MPDAGs in which the total effect of $A$ on $Y$ is identified, by considering all orientation configurations of the edges connected to $A$; see also \citet{liu2020collapsible} for a more efficient algorithm for a single treatment. However, this is often not minimal. For instance, to estimate the total effect of $A$ on $Y$ given MPDAG $\g$ in \cref{fig:example_IDA}(a), the IDA methods would enumerate four graphs listed in \cref{fig:example_IDA}(c) and thus report four estimates. But in fact, there are only three distinct total effects, which correspond to the MPDAGs listed in \cref{fig:example_IDA}(d). 

In this paper, we characterize the minimal additional edge orientations needed to identify a given total effect. Based on this characterization, we develop a recursive algorithm that outputs the minimal set of possible total effects along with the corresponding MPDAGs. Our results hold nonparametrically, that is, without assuming a particular type of data generating mechanism such as linearity. Furthermore, our results can be used in conjunction with recent developments on efficient effect estimators \citep{henckel,rotnitzky2019efficient,guo2020efficient} to produce a set of informative estimates.

\section{Preliminaries} \label{sec:prelim}

Throughout the paper we consider a random vector $X$, indexed by $V = \{V_1, \dots, V_p\}$, that is  $X = X_V$, such that each variable $X_{V_i}$ is represented by node $V_i$ in a graph $\g = (V,E,U)$.

\paragraph{Graphs, nodes and random variables.} A partially directed graph $\g= (V,E,U) $ consists of a set of nodes $ \vars=\left\lbrace V_{1},\dots,V_{p}\right\rbrace$ for $p \ge 1$, a set of directed ($\rightarrow$) edges $E $ and a set of undirected ($-$) edges $U$.

\paragraph{Induced subgraph.} An \textit{induced subgraph} $\g_{V'} =(V', E', U')$ of $\g= (V,E, U) $ consists of  $V' \subseteq V$, $E' \subseteq E$, and  $U' \subseteq U$ where $E'$ and $U'$ are all edges  between nodes in $V'$ that are in $E$ and $U$ respectively. 

\paragraph{Paths.} A \textit{path}  $p = \langle V_1, \dots , V_k \rangle$, $k > 1$ from $V_1 \in A$ to $V_k \in Y$ in $\g$ is a sequence of distinct nodes, such that $V_i$ and $V_{i+1}$, $i \in \{1, \dots , k-1\}$ are adjacent in $\g$.  A path of the form $V_1 - \dots - V_k$ is an undirected path and a path of the form $V_1 \rightarrow \dots \rightarrow V_k$ is a causal path. Additionally, $p$ is a \textit{possibly causal path} in $\g$ if no edge $V_{i} \leftarrow V_{j}, 0 \le i < j \le k$ is in $\g$. Otherwise, $p$ is a \textit{non-causal path} in $\g$ (see Definition 3.1 and Lemma 3.2 of \citealp{perkovic17}).
A path from node set $A$ to  $Y$ is \textit{proper} with respect to $A$ when only its first node is in $A$.  

\paragraph{Colliders, shields, and definite status paths.} If a path $p$ contains $V_i \rightarrow V_j \leftarrow V_k$ as a subpath, then $V_j$ is a \textit{collider} on $p$. A path $\langle V_{i},V_{j},V_{k} \rangle$ is an \emph{unshielded triple} if $ V_{i} $ and $ V_{k}$ are not adjacent. A path is \textit{unshielded} if all successive triples on the path are unshielded. A node $V_{j}$ is a \textit{definite non-collider} on a path $p$ if the  edge $V_i \leftarrow V_j$, or the edge $V_j \rightarrow V_k$ is on $p$, or if $V_{i} - V_j - V_k$ is a subpath of $p$ and $V_i$ is not adjacent to $V_k$. A node is of \textit{definite status} on a path if it is a collider, a definite non-collider or an endpoint on the path. A path $p$ is of definite status if every node on $p$ is of definite status.

\paragraph{d-connection, d-separation, and blocking.} A definite status path \textit{p} from node $A$ to node $Y$ is \textit{d-connecting} given a node set $Z$ ($A,Y \notin Z$) if every definite non-collider on $p$ is not in $Z$, and every collider on $p$ has a descendant in $Z$. Otherwise, $Z$ \textit{blocks} $p$.  If $Z$ blocks all definite status paths between $A$ and $Y$ in \mpdag{} $\g$, then $A$ is \textit{d-separated} from $Y$ given $Z$ in $\g$ and we write $A\perp_{\g} Y |Z$ \citep[Lemma C.1 of][]{henckel}. 

\paragraph{Probabilistic implications of d-separation.} Let $f$ be any observational density over $X$ consistent with an \mpdag{} $\g = (V,E,U)$. Let $A,Y$ and $Z$ be pairwise disjoint node sets in $V$. If $A$ and $Y$ are d-separated given $Z$ in $\g$, then $X_A$  and $X_Y$ are  conditionally independent  given $X_Z$ in the observational density $f$ \citep{lauritzen1990independence,pearl2009causality}.
Hence, all DAGs that encode the same d-separation relationships also encode the same conditional independences and are therefore \textit{Markov equivalent}.

\paragraph{Ancestral relationships.} If $A \to Y$ is in $\g$, then $A$ is a \textit{parent} of $Y$. If there is a causal path from node $A$ to node $Y$, then $A$ is an \textit{ancestor} of $Y$, and $Y$ is a \textit{descendant} of $A$.  If there is a possibly causal path from node $A$ to node $Y$, then $Y$ is a \textit{possible descendant} of $A$. We use the convention that every node is a descendant, ancestor, and  possible descendant of itself.
The sets of parents, ancestors, and  possible descendants of  a node $A$ in~$\g$ are denoted by $\Pa(A,\g)$, $\An(A,\g)$, $\PossDe(A,\g)$  respectively. For a set of nodes $A = \{A_1, \dots, A_k\}$, we let $\Pa(A,\g) =( \cup_{i =1}^{k}  \Pa(A_i,\g)) \setminus A$, $\An(A,\g) = \cup_{i=1}^{k}  \An(A_i,\g)$, and $\PossDe(A,\g) = \cup_{i=1}^{k} \PossDe(A_i,\g)$.

\paragraph{DAGs, PDAGs.}
A \textit{directed graph} contains only directed edges. 
A causal path from node $A$ to node $Y$ and  $Y\to A$ form a \textit{directed cycle}.  
A directed graph without directed cycles is a \textit{directed acyclic graph $(\DAG)$}.  A \textit{partially directed acyclic graph PDAG} is a partially directed graph without directed cycles.

\paragraph{Observational, interventional densities, and causal DAGs.} We consider \textit{do-interventions} $\Do(X_A =x_a)$ (for $A\subseteq V$), or $\Do(x_a)$ for shorthand, which represent outside interventions that set $X_a$ to a fixed value $x_a$. 
We call a density $f$ of $X$  under no intervention \textit{an observational density}.  
An observational density $f(x)$  is \textit{consistent} with a $\DAG$ $\g[D] =(V,E, \emptyset)$ if  $f(x)= \prod_{i=1}^{p} f(x_{v_i}|x_{\pa(v_{i},\g[D])})$ \citep{pearl2009causality}. 

A density $f(x|\Do(x_{a}))$  under intervention $\Do(X_A =x_a)$, $A \subseteq V$ is called an \textit{interventional density}.
An interventional density $f(x|\Do(x_{a}))$ is consistent with a DAG $\g[D] = (V,E,\emptyset)$ if there is an observational density $f$ consistent with  $\g[D]$ such that 
\begin{align}
f(x|\Do(x_{a})) =\prod_{\substack{i=1 \\ V_i \notin A}}^{p} f(x_{v_i}|x_{\pa(v_{i},\g[D])}),
\label{eq11}
\end{align}
for values $x_{\pa(v_i,\g[D])}$ of $X_{\Pa(V_i,\g[D])}$ that are consistent with $x_a$.
\cref{eq11} is known as the truncated factorization formula \citep{pearl2009causality}, manipulated density formula \citep{spirtes2000causation} or the g-formula \citep{robins1986new}. 

A DAG $\g[D] = (V,E,\emptyset)$ is \textit{causal} for a random vector $X$ if the observational and all interventional densities over $X$ are consistent with $\g[D]$.

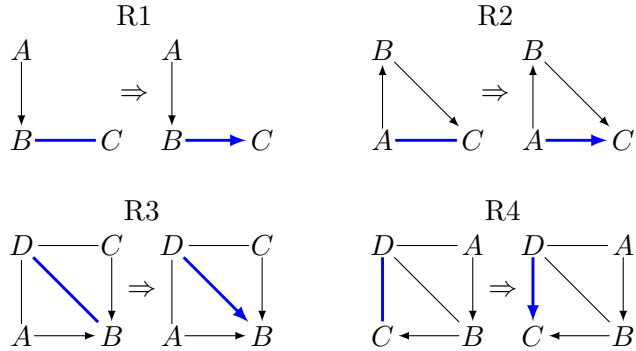
\begin{figure}[!tb]
\vspace{-.2cm}
\centering
\begin{tikzpicture}[->,>=latex,shorten >=1pt,auto,node distance=1.2cm,scale=1,transform shape]
  \tikzstyle{state}=[inner sep=0.5pt, minimum size=5pt]

\node[state] (Xia) at (0,0) {$B$};
  \node[state] (Xka) at (0,1.2) {$A$};
  \node[state] (Xja) at (1.2,0) {$C$};

  \path (Xka) edge (Xia);
 \draw[-,line width=1.1pt,blue]
          (Xia) edge (Xja);

  \coordinate [label=above:R1] (L1) at (1.5,1.4);
  \coordinate [label=above:$\Rightarrow$] (L2) at (1.5,0.4);

  \node[state] (Xib) at (2,0) {$B$};
  \node[state] (Xkb) at (2,1.2) {$A$};
  \node[state] (Xjb) at (3.2,0) {$C$};

  \path (Xkb) edge (Xib);
   \draw[->,line width=1.1pt,blue]   (Xib) edge (Xjb);

\node[state] (Xic) at (4.8,0) {$A$};
  \node[state] (Xkc) at (4.8,1.2) {$B$};
  \node[state] (Xjc) at (6,0) {$C$};

  \path (Xic) edge (Xkc)
          (Xkc) edge (Xjc);
   \draw[-,line width=1.1pt,blue]
          (Xic) edge (Xjc);

  \coordinate [label=above:R2] (L2) at (6.3,1.4);
  \coordinate [label=above:$\Rightarrow$] (L2) at (6.3,0.4);

  \node[state] (Xid)  at (6.8,0) {$A$};
  \node[state] (Xkd)  at (6.8,1.2) {$B$};
  \node[state] (Xjd)  at (8,0) {$C$};

  \path (Xid) edge (Xkd)
          (Xkd) edge (Xjd);
   \draw[->,line width=1.1pt,blue]  (Xid) edge (Xjd);

\node[state] (Xie) at (0,-1.4) {$D$};
  \node[state] (Xke) at (1.2,-1.4) {$C$};
  \node[state] (Xle) at (0,-2.6) {$A$};
  \node[state] (Xje) at (1.2,-2.6) {$B$};

  \path (Xke) edge (Xje)
          (Xle) edge (Xje);

  \draw[-,line width=1.1pt,blue]  (Xie) edge (Xje);
  \path[-]
          (Xke) edge (Xie)
          (Xle) edge (Xie);

  \coordinate [label=above:R3] (L3) at (1.6,-1.2);
    \coordinate [label=above:$\Rightarrow$] (L3) at (1.6,-2.2);

  \node[state] (Xif) at (2,-1.4) {$D$};
  \node[state] (Xkf) at (3.2,-1.4) {$C$};
  \node[state] (Xlf) at (2,-2.6) {$A$};
  \node[state] (Xjf) at (3.2,-2.6) {$B$};

  \draw[->,line width=1.1pt,blue]  (Xif) edge (Xjf);

  \path (Xkf) edge (Xjf)
          (Xlf) edge (Xjf);
  \path[-]
          (Xkf) edge (Xif)
          (Xlf) edge (Xif);

\node[state] (Xig) at (4.8,-1.4) {$D$};
  \node[state] (Xjg) at (6,-1.4) {$A$};
  \node[state] (Xkg) at (4.8,-2.6) {$C$};
  \node[state] (Xlg) at (6,-2.6) {$B$};

  \path (Xlg) edge (Xkg)
          (Xjg) edge (Xlg);
  \draw[line width=1.1pt,blue,-]
           (Xig) edge (Xkg);
  \path[-]
          (Xig) edge (Xjg)
          (Xig) edge (Xlg);

  \coordinate [label=above:R4] (L4) at (6.4,-1.2);
  \coordinate [label=above:$\Rightarrow$] (L4) at (6.4,-2.2);

  \node[state] (Xih) at (6.8,-1.4) {$D$};
  \node[state] (Xjh) at (8,-1.4) {$A$};
  \node[state] (Xkh) at (6.8,-2.6) {$C$};
  \node[state] (Xlh) at (8,-2.6) {$B$};

  \draw[blue,line width=1.1pt,->]  (Xih) edge (Xkh);
  \path (Xlh) edge (Xkh)
        (Xjh) edge (Xlh);
  \path[-]
          (Xih) edge (Xjh)
           (Xih) edge (Xlh);

\end{tikzpicture}
\caption{The orientation rules from~\cite{meek1995causal}. If the graph on the left-hand side of a rule is an induced subgraph of a PDAG $\g$, then \textit{orient} the blue undirected edge ({\color{blue} {\bfseries $-$}}) as shown on the right-hand side.} \label{fig:orientationRules}
\end{figure}

\paragraph{CPDAGs and MPDAGs.}   All DAGs that encode the same set of conditional independences  are \textit{Markov equivalent} and form a \textit{Markov equivalence class} of DAGs, which can be \textit{represented} by a \textit{completed partially directed acyclic graph} (CPDAG) \citep{meek1995causal,andersson1997characterization}.    A PDAG $\g$ is a \textit{maximally oriented} PDAG (\mpdag{}) if and only if the edge orientations in $\g$ are complete under rules R1-R4 in  Figure \ref{fig:orientationRules} \citep{meek1995causal}. An MPDAG is also known as CPDAG with background knowledge \citep{meek1995causal}. As such, both a DAG and a CPDAG can be seen as special cases of an \mpdag{}. Any graph in this paper can hence be labeled an MPDAG.

\paragraph{$\g$ and $[\g]$.} A DAG $\g[D] = (V,E,\emptyset)$ is \textit{represented} by \mpdag{} $\g = (V,E',U')$ if $\g[D]$ and $\g$ have the same adjacencies, same unshielded colliders and if $E' \subseteq E$ \citep{meek1995causal}.
If $\g$ is an \mpdag{}, then $[\g]$ denotes the set of all $\DAG$s represented by $\g$. 
An MPDAG $\g^{'}$ is said to be \textit{represented} by another MPDAG $\g$ if $[\g'] \subseteq [\g]$.

\paragraph{Causal MPDAGs.} An observational or interventional density is consistent with MPDAG $\g$ if it is consistent with a DAG $\g[D]$ in $[\g]$.
An MPDAG $\g = (V,E,U)$ is \textit{causal} if it represents the causal DAG.

\paragraph{Concatenation.} We denote the concatenation of paths by $\oplus$, so that for a path $p = \langle V_1,V_2,\dots,V_m \rangle$, $p = p(V_1, V_r) \oplus p(V_r, V_m)$, for $1\le r\le m$.

\paragraph{Buckets and bucket decomposition \citep{perkovic20}.}  A node set $A$, $A \subseteq V$  is an \textit{undirected connected set} in  $\g = ( V,E,U)$  if for every two distinct nodes $A_i, A_j \in A$, $A_i - \dots - A_j$ is in $\g$.  If node set $B$, $B \subseteq D \subseteq V$, is a maximal undirected connected subset of $D$ in $\g = (V,E,U)$, we call $B$ a \textit{bucket} in ${D}$. Additionally, $D$ can be partitioned into $D=D_1\,\cup \dots \cup\,D_K$, where each $D_k$, $k \in \{1, \dots , K\}$ is a bucket in $D$ and $D_i \cap D_j = \emptyset$ for $i \neq j$. We call the above partitioning of $D$ into buckets \textit{the bucket decomposition}.
Furthermore, $D_1, \dots, D_K$ can be ordered in such a way that if $D_1 \rightarrow D_2$ and $D_1 \in D_i$, $D_2 \in D_j$, then $i < j$; see PCO algorithm of \cite{perkovic20}.

\section{Main Results}

A total effect of $A$ on $Y$ is generally defined as some functional of the interventional density $f(x_y|\Do(X_A = x_a))$ (or $f(x_y |\Do(x_a))$ for short), such as $\dd \E [X_y | \Do(x_a)] / \dd x_a$ for continuous treatments and $\E [X_y | \Do(x_a=1)] - \E [X_y | \Do(x_a=0)]$ for binary treatments; see, e.g., \citet[Ch.\ 1]{hernanrobins}. 
For the common definitions, the total effect of $A$ on $Y$ is \textit{identified} in MPDAG $\g$ if and only if $f(x_y|\Do(x_a))$ can be identified from any observational density $f(x)$ consistent with $\g$ \citep{galles1995testing,perkovic20}. 
In this section, we show how to identify all MPDAGs represented by a given MPDAG $\g =(V,E,U)$ that have distinct identification maps for $f(x_Y|\Do(x_A))$, $A,Y \subseteq V$. Formally, the identification is a map from the space of observational densities that are consistent with $\g$ to the space of conditional kernels:
\[ f \in \mathcal{P}(\g) \mapsto f(x_Y | \Do(x_A)) \in \mathcal{K}(\mathcal{X}_A, \mathcal{X}_Y), \]
where $\mathcal{K}(\mathcal{X}_A, \mathcal{X}_Y)$ is the set of densities on the domain of $Y$ indexed by $A$. 
The identification (see \cref{thm:id-formula}) is possible if and only if $\g$ meets the following graphical condition. 

\begin{theorem}[{Identifiability condition of \citet{perkovic20}}] \label{thm:id-criterion}
Let $\g = (V,E,U)$ be a causal MPDAG for a random vector $X$. Further, let $A$ and $Y$ be disjoint node sets in $\g$. 
The total  effect of $A$ on $Y$ is identified in $\g$ if and only if every proper possibly causal path from $A$ to $Y$ starts with a directed edge in $\g$.
\end{theorem}

It follows that, if a total effect of $A$ on $Y$ is \emph{not} identified given MPDAG $\g$, then there is at least one proper possible causal path from $A$ to $Y$ in $\g$ that starts with an undirected edge $A_1 - V_1$ for  $A_1 \in A$ and $V_1 \in V$. 
Therefore, to identify the total effect, one can enumerate all the valid combinations of orientations just for the undirected edges of this type. 

As an example, consider MPDAG $\g$ in \cref{fig:example_IDA}(a). 
Paths $A - V_1 - Y$ and $A - V_2 - V_1 - Y$ are two proper possibly causal paths from $A$ to $Y$ in $\g$ that start with an undirected edge. 
There are four ways to orient the two starting edges: 
\begin{equation*}
\begin{split}
& R_1 = \{A \leftarrow V_2, A \rightarrow V_1\}, \quad R_2 = \{A \rightarrow V_2, A \rightarrow V_1\} \\
& R_3 = \{A \rightarrow V_2, A \leftarrow V_1\}, \quad R_4 =\{A \leftarrow V_2, A \leftarrow V_1\}.
\end{split}
\end{equation*}
Using algorithm \texttt{MPDAG}($\g$, $R_i$) for $i=1,\dots,4$ (\citealp{meek1995causal,perkovic17}, see Algorithm \ref{alg:mpdag} in the Appendix), which adds orientations $R_i$ and then completes the rules of \cite{meek1995causal} in \cref{fig:orientationRules}, we obtain the four MPDAGs listed in \cref{fig:example_IDA}(c). 
Now the effect of $A$ on $Y$ can be identified and estimated under each of the four MPDAGs.

This procedure already improves over the current standard IDA and joint-IDA algorithms \citep{maathuis2009estimating,nandy2017estimating,perkovic17,witte2020efficient} because orientations of fewer edges are considered --- IDA and joint-IDA orient \emph{all} undirected edges connected to $A$. 
However, we claim that it suffices to consider even fewer edges. 
The next theorem characterizes the minimal amount of edge orientation needed to identify a total effect. Its proof is left to the Appendix. 

\begin{theorem}\label{thm:min-bg}
Let $\g = (V,E,U)$ be a causal MPDAG. Let $A$ and $Y$ be disjoint node sets in $\g$ such that the total effect of $A$ on $Y$ is not identified given $\g$.
Suppose $p=\langle A_1, V_1, \dots, Y_1 \rangle$ for $A_1 \in A$, $Y_1 \in Y$ is a \emph{shortest} proper possibly causal path from $A$ to $Y$ such that $A_1 - V_1$. 
Then the total effect of $A$ on $Y$ is not identified in any MPDAG ${\g}^{*}$ that is represented by $\g$ and contains the undirected edge $A_1 -V_1$.
\end{theorem}

In the above, we say that ${\g}^{\ast}$ is represented by $\g$ if all DAGs represented by ${\g}^{\ast}$ are also represented by $\g$. 
Implicit in \cref{thm:min-bg} is the fact that when there are more than one paths that violate the identifiability condition (\cref{thm:id-criterion}), the order in which we orient them matters. In particular, the/a shortest path should be oriented first. 

Consider again $\g$ in \cref{fig:example_IDA}(a).
Because $A - V_1 - Y$ is shorter than $A - V_2 - V_1 - Y$, edge $A - V_1$ should be oriented first. 
In fact, as soon as this edge is oriented as $A \leftarrow V_1$, the acyclicity of the underlying DAG renders both paths $A \leftarrow V_1 - Y$ and $A - V_2 - V_1 -Y$ as non-causal from $A$ to $Y$; see the second row of \cref{fig:example_IDA}(d).

This characterization naturally leads to a recursive algorithm \texttt{IDGraphs} (Algorithm \ref{alg:computegraphs}). The \texttt{IDGraphs} algorithm takes MPDAG $\g$ and node sets $A$, $Y$ as input and outputs a finite set of MPDAGs $\{\g_1, \dots, \g_n\}$ that partition $\g$, such that (i) the total effect of $A$ on $Y$ is identified in every $\g_i$ and (ii) the effects identified from $\g_i$ and $\g_j$ are different for $i \neq j$. Property (ii), formally stated in \cref{thm:dif-effects}, shows that the enumeration is minimal. Additionally, by construction, it holds that $n \leq 2^{m(\g)}$, where $m(\g)$ is the number paths that violate the condition in \cref{thm:id-criterion}. 

To explain how the algorithm works, consider again the example in \cref{fig:example_IDA}(a). As we have already seen, \texttt{IDGraphs} first orients edge $A - V_1$ in $\g$. We obtain $\g_1 = \texttt{MPDAG}(\g, \{A_1 \rightarrow V_1\})$ and $\g_2 = \texttt{MPDAG}(\g, \{A_1 \leftarrow V_1\})$. Note that the effect is already identified in $\g_2$ (second row of \cref{fig:example_IDA}(d)) so $\g_2$ appears in the output. The effect in $\g_1$ is still not identified and the algorithm proceeds to orient edge $A - V_2$, which leads to the other two graphs in the output (first row of \cref{fig:example_IDA}(d)). 

\begin{algorithm}[tb]
 \Input{MPDAG $\g$, disjoint node sets $A$ and $Y$}
 \Output{the minimal set of MPDAGs with identified effects that partition $\g$ }
 \vspace{.2cm}
   \SetAlgoLined
 \uIf{$\g$ satisfies the condition in \cref{thm:id-criterion}}
 { \Return $\g$;\\
 } \Else{
 Let $A_1 -V_1 $ be an edge in $\g$ that satisfies Theorem \ref{thm:min-bg};\\
 ${\g}_{1} = $ MPDAG($\g$, $\{A_1 \to V_1\}$);\\
 ${\g}_{2} = $ MPDAG($\g$, $\{A_1 \leftarrow V_1\}$);\\
 \Return $\{\texttt{IDGraphs}(A,Y,{\g}_{1}), \texttt{IDGraphs}(A,Y,{\g}_{2})\}$;\\
 }
\caption{\texttt{IDGraphs}}
\label{alg:computegraphs}
\end{algorithm}

Note that, however, if $A - V_2$ was oriented before $A - V_1$, we would arrive at four MPDAGs (\cref{fig:example_IDA}(c)) instead of three!

The next theorem summarizes the theoretical guarantees of \texttt{IDGraphs}. We prove the minimality constructively; see the Appendix for details. 
\begin{theorem} \label{thm:dif-effects}
Suppose $\g$ is a causal MPDAG and $A,\, Y$ are two disjoint node sets in $\g$.
Let $L = \{\g_1, \dots, \g_n\}$ be the output of $\texttt{IDGraphs}(A, Y, \g)$. 
Then the following statements hold. 
\begin{enumerate}[label=(\roman*)]
\item The total effect of $A$ on $Y$ is identified in each $\g_i$.
\item For any $i \neq j$, there exists an observational density $f$ that is consistent with $\g$ such that the effect identified from $f$ in $\g_i$ is different from the effect identified from $f$ in $\g_j$. 
\item $L$ is a partition of $\g$ in terms of DAGs represented. 
\end{enumerate}
\end{theorem}
\begin{proof}[Proof sketch]
Here we sketch out the proof of (ii); see the Appendix for details. For two MPDAGs $\mathcal{G}_1$ and $\mathcal{G}_2$ output by \texttt{IDGraphs}, we \emph{construct} a density $f$ that factorizes according to $\mathcal{G}_1$ (and $\mathcal{G}_2$, due to Markov equivalence) but such that  $\E[X_Y | \text{do}(X_A = \bm{1})]$ has different values under $\mathcal{G}_1$ and $\mathcal{G}_2$. There are two steps in this construction. 

First, we establish some graphical differences between $\g_1$ and $\g_2$ that stem from the application of Thm 3.2 in the \texttt{IDGraphs} algorithm. Consider representing the recursion of \texttt{IDGraphs} as a binary tree and let $\mathcal{G}^{\ast}$ be the lowest common ancestor of $\mathcal{G}_1$ and $\mathcal{G}_2$. WLOG, suppose $A_1 - V_1 \in \mathcal{G}^{\ast}$ but $A_1 \rightarrow V_1 \in \mathcal{G}_1$, $A_1 \leftarrow V_1 \in \mathcal{G}_2$. Let $p^{\ast}$ be the shortest possibly causal path from $A$ to $Y$ in $\mathcal{G}^{\ast}$ that starts with $A_1 - V_1$. Then difference in terms of $p^{\ast}$ between $\mathcal{G}_1$ and $\mathcal{G}_2$ can be categorized into two cases, depending on if $p^{\ast}$ is shielded, given by \cref{lemma:auxillaryv1} in the Appendix. 
For each case, we \emph{parametrize} two linear Gaussian DAGs $\mathcal{D}_1 \in [\mathcal{G}_1]$ and $\mathcal{D}_2 \in [\mathcal{G}_2]$ such that their observed distributions are identical (by matching the first two moments) but values of $\E[X_{Y}| \Do(X_A = \mb{1})]$ are different. 
\end{proof}

\subsection{Examples} \label{sec:examples}
Suppose that a total effect of interest is not identified in an MPDAG $\g$. In order to obtain a set of possible total effects, one needs to enumerate the MPDAGs represented by $\g$ in which the total effect is identified. In below, we consider four methods that have appeared in our discussion so far, listed from the most computationally demanding to the least.

\begin{enumerate}[label = \cctext*,align=left]
\item\label{method1}  List all DAGs represented by $\g$. This is adopted by the global IDA algorithm of \citet{maathuis2009estimating}.

\item\label{method2}  List MPDAGs corresponding to all valid orientations for undirected edges attached to $A$. This is adopted by the local/semi-local IDA algorithms of \citet{maathuis2009estimating,perkovic17,witte2020efficient} and the joint-IDA algorithm of \citet{nandy2017estimating}.

\item\label{method3} List MPDAGs corresponding to  all valid combinations of edge orientations for edges $A_1 - V_1$, $A_1 \in A$, $V_1 \in V$, such that $V_1$ is on a proper possibly causal path from $A$ to $Y$ in $\g$.
\item\label{method4} Use \texttt{IDGraphs}($A,Y, \g$).
\end{enumerate}
These methods are compared through two examples, one for point intervention ($|A|=1$) and one for joint intervention ($|A|=2$).

\begin{example}[$|A| = 1$] \label{exam:point}
Consider again $\g$ in \cref{fig:example_IDA}(a). Recall that \ref{method1} lists 7 DAGs shown in \cref{fig:example_IDA}(b) and \ref{method4} lists 3 MPDAGs shown in \cref{fig:example_IDA}(d).
And as discussed earlier, both \ref{method2} and \ref{method3} would orient edges $A - V_2$ and $A - V_1$, and hence list 4 MPDAGs shown in \cref{fig:example_IDA}(d).
\end{example}

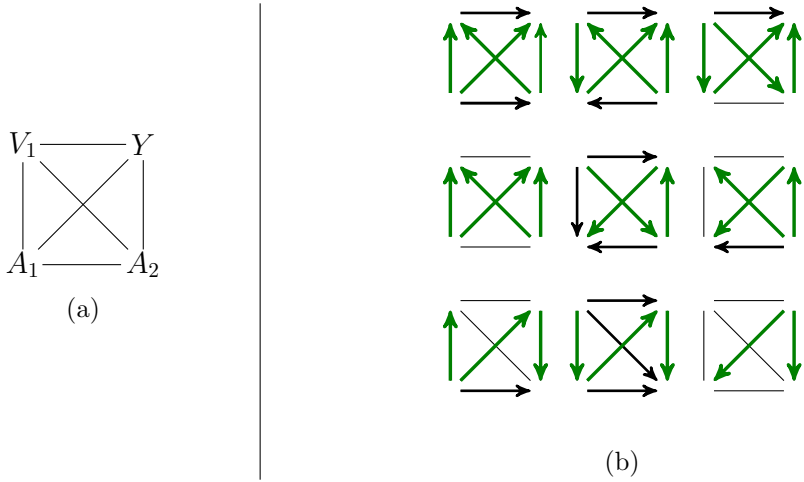
\begin{figure}[htb]
\tikzstyle{every edge}=[draw,>=stealth',->]
\newcommand\dagvariant[1]{\begin{tikzpicture}[xscale=.6,yscale=.6]
\node (a) at (0,0) {};
\node (d) at (0,2) {};
\node (b) at (2,0) {};
\node (c) at (2,2) {};
\draw (a) edge [-] (b);
\draw (b) edge [-] (c);
\draw (d) edge [-] (a);
\draw (d) edge [-] (c);
\draw (a) edge [-] (c);
\draw (b) edge [-] (d);
\draw #1;
\end{tikzpicture}}

\centering
\begin{subfigure}{.3\columnwidth}
  \centering
\begin{tikzpicture}[->,>=latex,shorten >=1pt,auto,node distance=0.8cm,scale=.8,transform shape]
  \tikzstyle{state}=[inner sep=1pt, minimum size=12pt]

\node[state] (Xia) at (0,0) {\Large $A_1$};
  \node[state] (Xka) at (0,2) {\Large $V_1$};
  \node[state] (Xja) at (2,0) {\Large  $A_2$};
  \node[state] (Xsa) at (2,2) {\Large  $Y$};

  \draw (Xka) edge [-] (Xsa);
  \draw (Xja) edge [-] (Xka);
\draw 		(Xia) edge  [-] (Xka);
 \draw   	(Xja) edge [-] (Xsa);
 \draw    	(Xia) edge [-] (Xja);
\draw    	(Xia) edge [-] (Xsa);
\end{tikzpicture}
\caption{}
  \label{mpdag2}
  \hspace{.2cm}
\end{subfigure}
\vrule
\begin{subfigure}{.65\columnwidth}
\begin{center}
\dagvariant{
(b)  edge [->,line width=1pt,color=ao]  (c)
(a) edge [->,line width=1.1pt] (b)
(b) edge [->,line width=1.3pt,color=ao] (d)
(a) edge [->,line width=1.3pt,color=ao] (c)
(a) edge [->,line width=1.3pt,color=ao] (d)
(d) edge [->,line width=1pt] (c)
}\hspace{.1cm}
\dagvariant{
(b)  edge [->,line width=1.3pt,color=ao]  (c)
(b) edge [->,line width=1pt] (a)
(a) edge [->,line width=1.3pt,color=ao] (c)
(d) edge [->,line width=1.3pt,color=ao] (a)
(d) edge [->,line width=1pt] (c)
(b) edge [->,line width=1.3pt,color=ao] (d)
}\hspace{.1cm}
\dagvariant{
(b)  edge [->,line width=1.3pt,color=ao]  (c)
(a) edge [->,line width=1.3pt,color=ao] (c)
(d) edge [->,line width=1.3pt,color=ao] (a)
(d) edge [->,line width=1.3pt,color=ao] (b)
(d) edge [->,line width=1pt] (c)
}\\
\end{center}
\begin{center}
\dagvariant{
(b)  edge [->,line width=1.3pt,color=ao]  (c)
(a) edge [->,line width=1.3pt,color=ao] (c)
(b) edge [->,line width=1.3pt,color=ao] (d)
(a) edge [->,line width=1.3pt,color=ao] (d)
}\hspace{.1cm}
\dagvariant{
(b)  edge [->,line width=1pt]  (a)
(b)  edge [->,line width=1.3pt,color=ao]  (c)
(c) edge [->,line width=1.3pt,color=ao] (a)
(d) edge [->,line width=1.3pt,color=ao] (b)
(d) edge [->,line width=1pt] (a)
(d) edge [->,line width=1pt] (c)
}\hspace{.1cm}
\dagvariant{
(b)  edge [->,line width=1pt]  (a)
(b)  edge [->,line width=1.3pt,color=ao]  (c)
(c) edge [->,line width=1.3pt,color=ao] (a)
(b) edge [->,line width=1.3pt,color=ao] (d)
}\\
\end{center}
\begin{center}
\dagvariant{
(a)  edge [->,line width=1pt]  (b)
(c)  edge [->,line width=1.3pt,color=ao]  (b)
(a)  edge [->,line width=1.3pt,color=ao]  (c)
(a) edge [->,line width=1.3pt,color=ao] (d)
}\hspace{.1cm}
\dagvariant{
(a) edge [->,line width=1pt] (b)
(d) edge [->,line width=1pt] (c)
(d) edge [->,line width=1pt] (b)
(c)  edge [->,line width=1.3pt,color=ao]  (b)
(a)  edge [->,line width=1.3pt,color=ao]  (c)
(d) edge [->,line width=1.3pt,color=ao] (a)
}\hspace{.1cm}
\dagvariant{
(c)  edge [->,line width=1.3pt,color=ao]  (a)
(c)  edge [->,line width=1.3pt,color=ao]  (b)
}\\
\end{center}
\caption{}
\end{subfigure}
\caption{\cref{exam:joint}. (a) $\MPDAG$ $\g$, (b) all $\MPDAG$s with distinct causal identification formulas for $f(x_Y|\Do(x_{A_1}, x_{A_2}))$ represented by $\g$.} 
\label{fig:example_jointIDA}
\end{figure}

\begin{example}[$|A|=2$] \label{exam:joint}
Consider $\MPDAG$ $\g$ in \cref{fig:example_jointIDA}(a) and let $A = \{A_1, A_2\}$. \ref{method4} yields 9 MPDAGs listed in \cref{fig:example_jointIDA}(b), which correspond to 9 distinct possible effects of $(A_1, A_2)$ on $Y$. 
\ref{method2} would consider all valid combinations of orientations for edges $A_1 - A_2$, $A_1 -V_1$, $A_1 - Y$, $A_2 - V_1$, and $A_2 -Y$, resulting in 18 MPDAGs.
Lastly, \ref{method3} would consider all valid combinations of orientations for edges $A_1 - V_1$, $A_1 - Y$, $A_2 - V_1$, and $A_2 -Y$, resulting in 12 MPDAGs.
\end{example}
\subsection{Computational Complexity}
The computational complexities of various algorithms are summarized in \cref{tab:complexity}. Here $l(\mathcal{G})$ is the number of undirected edges incident to $A$. For the running time of collapsible IDA \citep{liu2020collapsible}, $O(|V| +|E|)$ reflects finding the neighbors of $A$ that are on a possible causal path to $Y$ and $r(\mathcal{G})$ is  the size of that subset. Both local IDA and collapsible IDA are only applicable when $|A|=1$. 

\begin{table}[!htb]
\caption{Computational complexity}
\label{tab:complexity}
\centering
\begin{tabular}{@{}lll@{}}
\toprule
local IDA & $O(2^{l(\mathcal{G})})$  \\
collapsible IDA &  $O((|V| + |E|)2^{r(\mathcal{G})})$ \\
semi-local/joint/optimal IDA & $O(2^{l(\mathcal{G})}) \text{poly}(|V|)$  \\
\texttt{IDGraphs} & $O(2^{m(\mathcal{G})}) \text{poly}(|V|)$ \\ \bottomrule
\end{tabular}
\end{table}

For more general settings, our \texttt{IDGraphs} algorithm is asymptotically on par with semi-local \citep{maathuis2009estimating,perkovic17}, joint \citep{nandy2017estimating}, and optimal \citep{witte2020efficient} variants of the IDA algorothm. The running time of these methods is bounded by $O(2^{l(\mathcal{G})}) \text{poly}(|V|)$, where $\text{poly}(|V|)$ time is used to complete the orientation rules of \citet{meek1995causal}. Similarly, the complexity of \texttt{IDGraphs} is $ O(2^{m(\mathcal{G})}) \text{poly}(|V|)$,
where $m(\mathcal{G})$ is the number of undirected edges incident to $A$ on a proper possibly causal path from $A$ to $Y$ (\cref{thm:id-criterion}). Clearly, we have $m(\mathcal{G}) \leq l(\mathcal{G})$. The number of recursions is bounded by $2^{m(\mathcal{G})}$ and the time for each recursion by $\text{poly}(|V|)$, which includes completing the orientation rules, checking the condition of \cref{thm:id-criterion}, and identifying the shortest path (\cref{thm:min-bg}) if the condition is not met. 

\begin{figure}[!htb]
\centering
\includegraphics[width=0.7\textwidth]{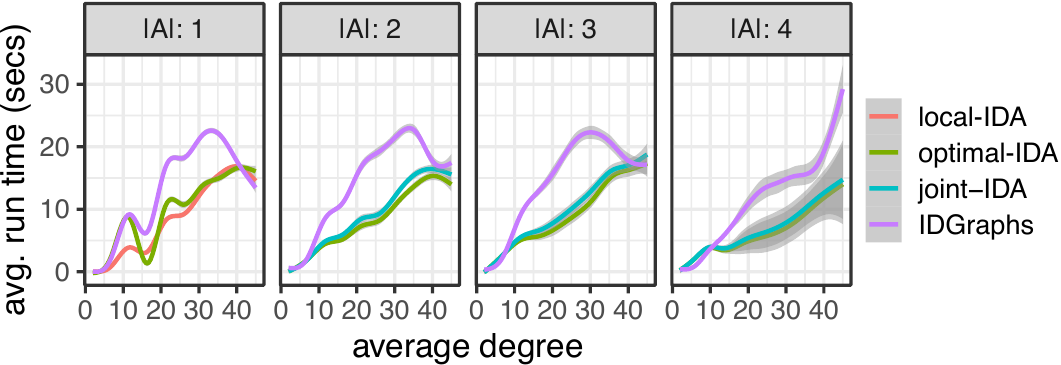}
\caption{Average run time versus the average node degree of the graph (shade: 95\% confidence band).}
\label{fig:runtime}
\end{figure}

In practice, from the simulations in \cref{sec:numerical}, we find that \texttt{IDGraphs} roughly costs twice the time of IDA type algorithms; see \cref{fig:runtime}. 

\subsection{Corollaries of the Main Result}
Having obtained a set of MPDAGs from the \texttt{IDGraphs} algorithm, the total effect can be estimated for each graph. Some recently developed efficient estimators can be employed, including the semiparametric efficient estimator of \citet{rotnitzky2019efficient} and the efficient least-squares estimator of \citet{guo2020efficient} under linearity assumptions. This strategy applies to any other causal quantity that is a functional of the interventional density.

\begin{theorem}[{Causal identification formula}, \citealp{perkovic20}]\label{thm:id-formula}
Let $\g$ be a causal MPDAG. Suppose $\g$ satisfies the identifiability condition in terms of the total effect of $A$ on $Y$ (Theorem \ref{thm:id-criterion}). 
Further, let $B =\An(Y, \g_{V \setminus A}) \setminus Y$ and let $B_1, \dots, B_k$ be the bucket decomposition of $B \cup Y$.
Then for any observational density $f$ that is consistent with $\g$, we have
\begin{equation} \label{eqs:g-formula}
f(x_y|\Do(x_a)) = \int \prod_{i=1}^{k} f(x_{b_i}| x_{\pa(b_i,\g)}) \dd x_b,
\end{equation}
where values of $x_{\pa(b_i,\g)}$ are consistent with $x_a$.
\end{theorem}

Let $L$ be the output of $\texttt{IDGraphs}(A, Y, \g)$. 
\begin{corollary} \label{corollary:causal-id-form}
Then there are no two graphs in $L$ that share the same formula \cref{eqs:g-formula}.
\end{corollary}

Another prominent method for identifying the interventional density is through covariate adjustment. 
The IDA algorithms of \citet{maathuis2009estimating, perkovic17, witte2020efficient} are based on covariate adjustment for causal linear models.  
See \citet{perkovic20} for the generalized adjustment criterion, which is necessary and sufficient  for covariate adjustment in MPDAGs and generalizes the well-known back-door formula of \citet{pearl1993bayesian}.

\begin{corollary} \label{cor:adjust}
Then there are no two MPDAGs in $L$ that share the same adjustment set relative to $(A,Y)$. Further, if $|A| =|Y|=1$, then there exists an adjustment set relative to $(A,Y)$ for each MPDAG in $L$.
\end{corollary}

\section{NUMERICAL RESULTS} \label{sec:numerical}
\begin{table*}[!ht]
\caption{Estimates of possible effects for the example in \cref{fig:simu-example}. Symbol $(a)^b$ denotes value $a$ with multiplicity $b$ in the multiset. Note that IDA (local) and joint-IDA return more values than the number of possible effects.}
\label{tab:simu-example}
\begin{center} \small
\begin{tabular}{llllllllll}
\toprule[1.2pt]
& $A_1$ on $Y$ (\cref{fig:simu-example}(c)) & & $A_1, A_2$ on $Y$ (\cref{fig:simu-example}(d))\\ 
\cline{2-2} \cline{4-4} 
true effect & 3 & & (2,1) \\
true possible effects & $\{3,\,2,\,1.8,\,0\}$ & & $\{(2,1),\,(3,0),\,(0,2),\,(0,0)\}$  \\
our method & $\{2.9,\,2.1,\,1.9,\,0\}$ & & $\{(2.1, 0.9),\,(2.9, 0),\,(0, 1.9),\,(0,0)\}$ \\ 
IDA (optimal) \citep{witte2020efficient} & $\{2.9,\,(2.1)^2,\,1.9,\,0\}$ & & $\{(2.1, 0.9)^6,\, (0,0)^2,\,(\texttt{NA}, \texttt{NA})^2\}$ \\
IDA (local) \citep{maathuis2009estimating} & $\{2.9,\,2.1,\,2.2,\,1.9,\,0\}$ & & --- \\
joint-IDA \citep{nandy2017estimating} & --- & & $\{(2.1, 0.9)^2,\,(2.2, 0.9),\,(1.9,1.1),\,(2.2, 1.1)^2,$ \\
& & & $ \quad (0, 1.9),\,(2.9, 0),\,(0,0)^2 \}$ \\
\bottomrule[1.2pt]
\end{tabular}
\end{center}
\end{table*}

\begin{figure*}[t]
\centering
\includegraphics[width=1.0\textwidth]{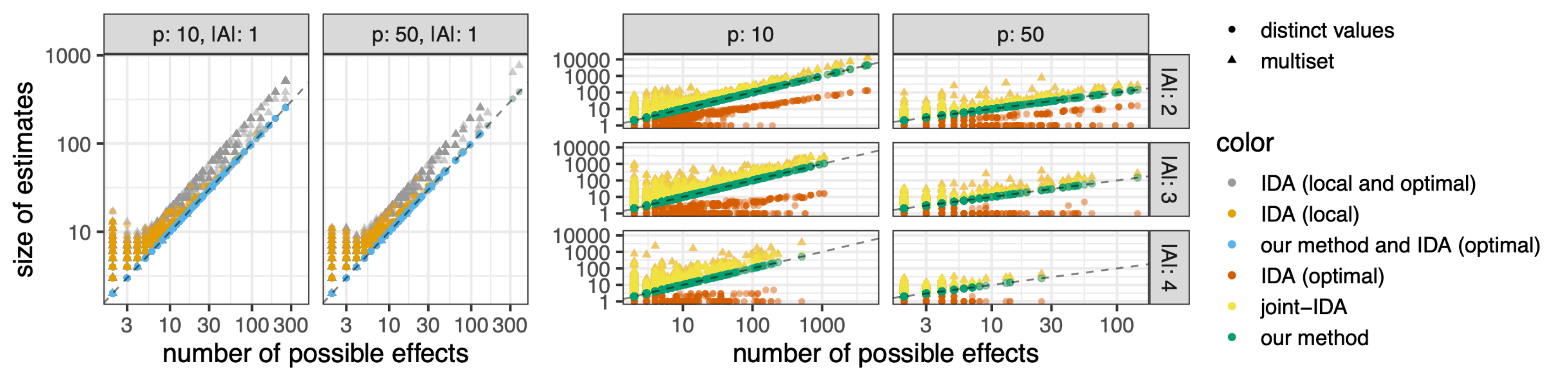}
\caption{Size of estimates vs. the number of possible effects on random instances (left: $|A|=1$, right: $|A|>1$). Both axes are in logarithmic scale. A dot on the graph represents either the size of the multiset ($\blacktriangle$) or the number of distinct values ($\bullet$).}
\label{fig:simu}
\end{figure*}

\begin{figure}[!hbt]
\tikzstyle{every edge}=[draw,>=stealth',->]
\newcommand\dagvariant[1]{\begin{tikzpicture}[rv/.style={circle, draw, minimum size=1.5mm, inner sep=0mm}, xscale=.6, yscale=0.6]
\node[rv] (a) at (1,2) {};
\node[rv] (b) at (1,1) {};
\node[rv] (v) at (0,0) {};
\node[rv] (y) at (2,0) {};
\draw #1;
\end{tikzpicture}}

\tikzstyle{every edge}=[draw,>=stealth',->]
\newcommand\dagempty[1]{\begin{tikzpicture}[xscale=.5,yscale=0.5]
\node (a) at (1,2) {};
\node (b) at (1,1) {};
\node (v) at (0,0) {};
\node (y) at (2,0) {};
\draw #1;
\end{tikzpicture}}

\centering
\begin{subfigure}{.4\columnwidth}
  \centering
\begin{tikzpicture}[->,>=latex,shorten >=1pt,auto,node distance=0.8cm,scale=1,transform shape]
  \tikzstyle{state}=[inner sep=1pt, minimum size=12pt]

\node[state] (a) at (1,2) {\large $A_1$};
  \node[state] (b) at (1,1) {\large $A_2$};
  \node[state] (v) at (0,0) {\large $V$};
  \node[state] (y) at (2,0) {\large $Y$};

\draw (a) edge [->, thick] (b);
\draw (a) edge [->, thick] (v);
\draw (a) edge [->, thick] (y);
\draw (b) edge [->, thick] (v);
\draw (b) edge [->, thick] (y);
\end{tikzpicture}
\caption{}
\end{subfigure}
\vrule
\hspace{0.2cm}
\begin{subfigure}{.4\columnwidth}
  \centering
\begin{tikzpicture}[->,>=latex,shorten >=1pt,auto,node distance=0.8cm,scale=1,transform shape]
  \tikzstyle{state}=[inner sep=1pt, minimum size=12pt]

\node[state] (a) at (1,2) {\large $A_1$};
  \node[state] (b) at (1,1) {\large $A_2$};
  \node[state] (v) at (0,0) {\large $V$};
  \node[state] (y) at (2,0) {\large $Y$};

\draw (a) edge [-, thick] (b);
\draw (a) edge [-, thick] (v);
\draw (a) edge [-, thick] (y);
\draw (b) edge [-, thick] (v);
\draw (b) edge [-, thick] (y);
\end{tikzpicture}
\caption{}
\end{subfigure}
\begin{subfigure}{.4\columnwidth}
\begin{center}
\dagvariant{
(a) edge [->,thick,color=ao] (b)
(a) edge [->,thick,color=ao] (v)
(a) edge [->,thick,color=ao] (y)
(b) edge [-,thick] (v)
(b) edge [-,thick] (y)
}
\dagvariant{
(a) edge [<-,thick,color=ao] (b)
(a) edge [-,thick] (v)
(a) edge [->,thick,color=ao] (y)
(b) edge [-,thick] (v)
(b) edge [->,thick,color=ao] (y)
}
\end{center}
\begin{center}
\dagvariant{
(a) edge [->,thick,color=ao] (b)
(a) edge [<-,thick,color=ao] (v)
(a) edge [->,thick,color=ao] (y)
(b) edge [<-,thick,color=ao] (v)
(b) edge [->,thick,color=ao] (y)
}
\dagvariant{
(a) edge [-,thick] (b)
(a) edge [->,thick,color=ao] (v)
(a) edge [<-,thick,color=ao] (y)
(b) edge [->,thick,color=ao] (v)
(b) edge [-,thick] (y)
}
\end{center}
\caption{}
\end{subfigure}
\vrule
\hspace{0.2cm}
\begin{subfigure}{.4\columnwidth}
\begin{center}
\dagvariant{
(a) edge [-,thick] (b)
(a) edge [-,thick] (v)
(a) edge [->,thick,color=ao] (y)
(b) edge [-,thick] (v)
(b) edge [->,thick,color=ao] (y)
}
\dagvariant{
(a) edge [->,thick,color=ao] (b)
(a) edge [->,thick,color=ao] (v)
(a) edge [->,thick,color=ao] (y)
(b) edge [->,thick,color=ao] (v)
(b) edge [<-,thick,color=ao] (y)
}
\end{center}
\begin{center}
\dagvariant{
(a) edge [<-,thick,color=ao] (b)
(a) edge [->,thick,color=ao] (v)
(a) edge [<-,thick,color=ao] (y)
(b) edge [->,thick,color=ao] (v)
(b) edge [->,thick,color=ao] (y)
}
\dagvariant{
(a) edge [-,thick] (b)
(a) edge [->,thick,color=ao] (v)
(a) edge [<-,thick,color=ao] (y)
(b) edge [->,thick,color=ao] (v)
(b) edge [<-,thick,color=ao] (y)
}
\end{center}
\caption{}
\end{subfigure}
\caption{Estimating possible total effects under a linear causal model: (a) the underlying causal DAG $\mathcal{D}$, (b) CPDAG $\mathcal{G}$ that represents $\mathcal{D}$, (c) $\texttt{IDGraphs}(A_1, Y, \mathcal{G})$, (d) $\texttt{IDGraphs}(\{A_1, A_2\}, Y, \mathcal{G})$.} 
\label{fig:simu-example}
\end{figure}
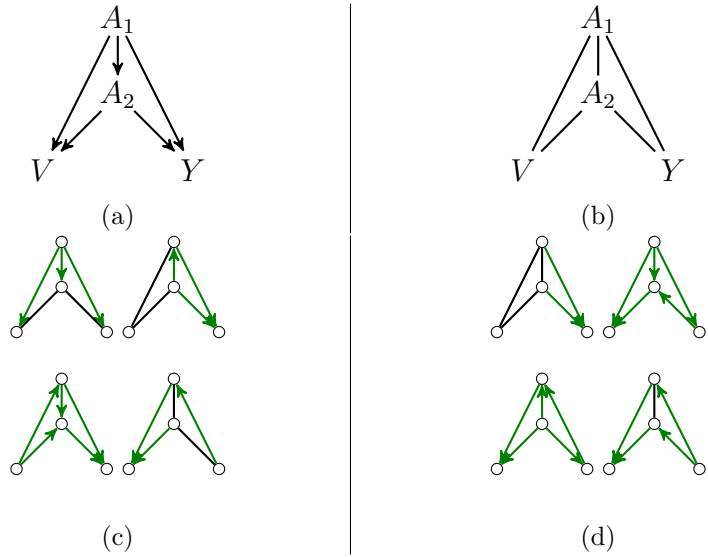

We present numerical results on estimating possible total effects under a linear causal model \citep{bollen1989structural}. To fix the idea, consider the example in \cref{fig:simu-example}. Suppose data is generated from an underlying linear model 
\begin{equation*}
\begin{split}
&X_{A_1} = \varepsilon_1, \quad X_{A_2} = \gamma_{12} X_{A_1} + \varepsilon_2, \\
&X_V = \gamma_{13} X_{A_1} + \gamma_{23} X_{A_2} + \varepsilon_3, \\
&X_Y = \gamma_{14} X_{A_1} + \gamma_{24} X_{A_2} + \varepsilon_4,
\end{split}
\end{equation*}
associated with causal DAG $\mathcal{D}$ given by (a). In this case, we set $\gamma_{14} = \gamma_{23} = 2$ and $\gamma_{12}=\gamma_{13}=\gamma_{24}=1$. The errors $\varepsilon_i$ for $i=1,\dots,4$ are drawn independently from the standard normal distribution. Suppose the causal DAG is known up to its CPDAG $\mathcal{G}$ (no added background knowledge), which is shown in (b). We consider estimating the total effect of $A_1$ on $Y$ (point intervention), and the total effect of $(A_1, A_2)$ on $Y$ (joint intervention). 

\cref{tab:simu-example} shows the estimates from 100 samples, where ``our method'' refers to applying the efficient estimator of \citet{guo2020efficient} to each graph returned by \texttt{IDGraphs}. In general, the IDA algorithms enumerate possible graphs where the effect is identified and return the estimates as a multiset; see \ref{method1} and \ref{method2} in \cref{sec:examples}. Further, the distinct values of the multiset can be taken as the estimates of possible effects. However, as we can see, one possible effect can correspond to more than one distinct values due to sampling variability. Moreover, \texttt{NA}'s are produced when applying IDA (optimal) to joint interventions due to nonexistence of valid adjustment sets \citep{perkovic18}. 

To examine these issues in more generality, we simulate random instances and compare the size of estimates to the true number of possible effects. Causal DAG $\mathcal{D}$ is generated by sampling from the Erd{\H{o}}s-R{\'e}nyi model and assigning a random causal ordering. We consider graphs of size $p=10$ and $p=50$, where the average degree $k$ is drawn from $\{2,\dots,8\}$ for the former and $\{2,\dots, 45\}$ for the latter. We take $\mathcal{G}$ to be the CPDAG of $\mathcal{D}$. Treatment variables $A$ and outcome $Y$ are randomly selected such that the total effect of $A$ on $Y$ is unidentified from $\mathcal{G}$. The size of $A$ varies from 1 to 4. For each instance, given its $\mathcal{G}$ and 500 independent samples generated by a corresponding linear causal model (with random coefficients and errors drawn from the standard normal), the possible effects of $A$ on $Y$ are estimated. 

The result is summarized in \cref{fig:simu} from roughly 55,000 random instances. For point interventions (left panel), IDA (local) produces duplicates, i.e., more distinct values than the actual number of possible effects, especially when the number of possible effects is small, whereas IDA (optimal) returns the correct number of distinct values. For joint interventions (right panel), the joint-IDA algorithm often suffers from an excessive amount of duplicates, while IDA (optimal), on the other hand, severely underreports the size of possible effects due to too many \texttt{NA}'s it produced --- note the logarithmic scale of both axes. Therefore, our \texttt{IDGraphs} algorithm, in conjunction with a statistically efficient estimator of an identified effect, should be used in place of IDA algorithms in both cases, to avoid unnecessary computational overhead and deliver causally informative estimates.

\section{DISCUSSION}
We have studied the set-identification of a total effect given that the underlying causal DAG is known up to a Markov equivalence class or its refinement. Existing enumerative approaches to this problem are often not minimal, which cause unnecessary computational overheads and undesirable statistical duplicates. To ensure minimality, there are two key ingredients. The first is to determine the class of DAGs that should be grouped together, which is given by the identifiability condition of \cref{thm:id-criterion}. This condition locates the set of ``problematic'' undirected edges that must be oriented, while the other undirected edges can be left intact. The second key ingredient is to determine the order in which the ``problematic'' edges are oriented. This order dependency is tricky because orienting one edge may imply the orientation of another edge due to the orientation rules \citep{meek1995causal}. Perhaps surprisingly, the optimal order can be determined, which is to orient the shortest ``problematic'' path first (\cref{thm:min-bg}). This naturally leads to  \texttt{IDGraphs}, a simple recursive algorithm that guarantees minimal enumeration (\cref{thm:dif-effects}). \texttt{IDGraphs} can be readily used in conjunction with recent developments in efficient estimation \citep{henckel,rotnitzky2019efficient,guo2020efficient} to deliver informative estimates of the true causal effect. From this perspective, our result can be viewed as \emph{separating} two sources of uncertainty --- identification and estimation~---~the two crucial steps in causal inference.

One may wonder whether this approach can be extended to allow latent variables. The latent variable IDA (LV-IDA) algorithm of \citet{malinsky2017estimating} employs a strategy similar to the IDA algorithm, but instead given a partial ancestral graph (PAG). A PAG represents a Markov equivalence class of maximal ancestral graphs (MAGs), which are obtained from DAGs by marginalizing over latent variables. 
For one obstacle in this setting, it still unclear how to incorporate background knowledge of edge orientations into a PAG. 
In addition, \citet{pmlr-v97-jaber19a} recently showed that if an effect is not identified given a PAG, then there is at least one MAG in the Markov equivalence class in which the effect is still not identified (see their Theorem 4) --- the same enumeration strategy will no longer work.

We conclude with a final remark. Strictly speaking, the local versions of IDA and joint-IDA algorithms only require the neighborhood information of $A$ instead of the whole MPDAG. Yet, it is unclear whether there are many situations where one knows the neighborhood without knowing more structures. 

\paragraph{Acknowledgment} FRG acknowledges the support from ONR Grant N000141912446.

\appendix
\counterwithin{figure}{section}
\counterwithin{theorem}{section}

\section{Additional Preliminaries} \label{appendix:prelim}

\begin{lemma}[Rules of the do-calculus, \citealp{pearl2009causality}]
\label{lem:do-rules}
Let $A,Y,Z$ and $W$ be pairwise disjoint (possibly empty) node sets in causal DAG $\g[D] = (V,E,\emptyset)$.
Let $\g[D]_{\overline{A}}$ denote the graph obtained by deleting all edges into $A$ from $\g[D]$. Similarly, let $\g[D]_{\underline{A}}$ denote the graph obtained by deleting all edges out of $A$ in~$\g[D]$ and let $\g[D]_{\overline{A}\underline{Z}}$ denote the graph obtained by deleting all edges into $A$ and all edges out of $Z$ in~$\g[D]$.

\textbf{Rule 1.} If $Y \perp_{\g[D]_{\overline{A}}} Z | A \cup W$, then 
$f(x_y | \Do(x_a), x_w) = f(x_y | \Do(x_a),x_z,x_w). $

\textbf{Rule 2.}  If $Y  \perp_{\g[D]_{\overline{A}\underline{Z}}} Z | A \cup W$, then
$ f(x_y | \Do(x_a),\Do(x_z), x_w) =f(x_y | \Do(x_a), x_z, x_w).$

\textbf{Rule 3.} If $Y \perp_{\g[D]_{\overline{{AZ(W)}}}}  Z | A \cup W$, then
$ f(x_y | \Do(x_a), x_w) =f(x_y | \Do(x_a), x_z, x_w), $
where $Z(W) =  Z \setminus \An(W, \g[D]_{\overline{A}})$.
\end{lemma}

\begin{lemma}[Lemma 3.6 of \citealp{perkovic17}]
 Let $A$ and $Y$ be distinct nodes in a \mpdag{} $\g$. If $p$ is a possibly causal path from $A$ to $Y$ in $\g$, then a subsequence $\pstar$ of $p$ forms a possibly causal unshielded path from $A$ to $Y$ in~$\g$.
\label{lemma:unshielded-analog-ema-uai}
\end{lemma}

\begin{lemma}[Wright's rule, \citealp{wright1921correlation}]
Let $X= AX + {\epsilon}$, where $A \in \mathbb{R}^{p \times  p}$, $X= X_V$, $|V| = p$,  and ${\epsilon} = (\epsilon_1,\dots, \epsilon_p)^T$ is a vector of mutually independent errors with means zero and proper variance such that  $\Var(X_i) = 1$, for all $i \in \{1, \dots , p\}$.
Let $\g[D] = (V,E,\emptyset)$, be the corresponding $\DAG$ For two distinct nodes $i, j \in V$, let $p_1, \dots, p_k$ be all paths between $i$ and $j$ in $\g[D]$ that do not contain a collider. Then $\Cov(X_i,X_j) = \sum_{r=1}^{k}\pi_r$, where $\pi_r$ is the product of all edge coefficients along path $p_r$, $r \in \{1,\dots, k\}$.
\label{lem:wright}
\end{lemma}

\begin{lemma} (See, e.g., \citealp[Theorem 3.2.4]{mardia1980multivariate})
Let $X = ({X_1},{X_2})$ be a $p$-dimensional multivariate Gaussian random vector with mean vector ${\mu} = ({\mu_1},\mu_2)$ and covariance matrix ${\Sigma} = \begin{bmatrix}
{\Sigma_{11}} & {\Sigma}_{12} \\
{\Sigma_{21}} & {\Sigma}_{22}
\end{bmatrix}$, so that $X_1$ is a $q$-dimensional multivariate Gaussian random vector with mean vector $\mu_1$ and covariance matrix ${\Sigma}_{11}$ and  ${X_2}$ is a $(p-q)$-dimensional multivariate Gaussian random vector with mean vector ${\mu_2}$ and covariance matrix $\Sigma_{22}$.
Then $\E[{X_2} |{X_1 = x_1}] = \mu_2+ {\Sigma}_{21}{\Sigma}_{11}^{-1} ({x_1} - {\mu_1})$.
\label{lem:mardia-condexp}
\end{lemma}

\begin{algorithm}[t]
 \Input{MPDAG $\g$, set of background knowledge edge orientations $R$.}
 \Output{MPDAG $\g'$ or FAIL.}
 \vspace{.2cm}
   \SetAlgoLined
Let $\g' = \g$\;
 \While{$R \neq \emptyset$}{
  Choose an edge $\{U \rightarrow V\}$ in $R$\;
  $R = R\setminus \{U \rightarrow V\}$\;
  \eIf{$\{U - V\}$ or $\{U \rightarrow V\}$ is in $\g'$}{
   Orient $\{U \rightarrow V\}$ in $\g'$\;
   Iterate the rules in \cref{fig:orientationRules} of until no more can be applied;
   }{
   FAIL\;
  }
}
\Return  $\g'$\;
\caption{\texttt{MPDAG} (see also \citealp{meek1995causal} and Algorithm 1 of \citealp{perkovic17})}
\label{alg:mpdag}
\end{algorithm}

\section{Proofs Of Main Results} \label{sec:proofs}
\begin{proofof}[\cref{thm:min-bg}]
Let $p = \langle A_1,  V_1, \dots, V_k = Y_1 \rangle$, $k \ge 1$, $A_1 \in A$, $Y_1 \in Y$.
If $k =1$, that is, $A_1 - Y_1$ is in $\g$, the proposition clearly holds. Hence, we will assume $k>1$.
Suppose for a contradiction that there is an $\MPDAG$ ${\g}^{*}$ represented by $\g$ such that $A_1-V_1$ is in  ${\g}^{*}$ and that the total effect of $A$ on $Y$ is identified in ${\g}^{*}$.  Further, let $\pstar$ be the path in ${\g}^{*}$ that corresponds to path $p$ in $\g$, so that $p$ and $\pstar$ are both sequences of nodes $\langle A_1, V_1, \dots, V_k = Y_1 \rangle$, $k >1$.

Since the total effect of $A$ on $Y$ is identified in ${\g}^{*}$,  and because $\pstar$ is a proper path from $A$ to $Y$ that starts with an undirected edge in ${\g}^{*}$, by \cref{thm:id-criterion}, $\pstar$ must be a non-causal path from $A_1$ to $Y_1$ in ${\g}^{*}$. We show that this implies that $A_1 - V_1 \leftarrow V_2$ and $A_1 \rightarrow V_2$ are in ${\g}^{*}$, which contradicts that ${\g}^{*}$ is an MPDAG (because orientations in ${\g}^{*}$ are not complete with respect to R2 in \cref{fig:orientationRules}).

We first show that any existing edge between $A_1$ and $V_i$, $i \in \{2, \dots , k\}$ in $\g$ is of the form $A_1 \rightarrow V_i$. Suppose that there is an edge between $A_1$ and $V_i$, in $\g$. This edge is cannot be of the form $A_1 \leftarrow V_i$, since that would imply that $p$ is a non-causal path in $\g$. This edge also cannot be of the form $A_1 - V_i$, because otherwise we can concatenate $A_1 - V_i$ and $p(V_i, Y_1)$ to construct a proper possibly causal path from $A$ to $Y$ in $\g$ that is shorter than $p$. Hence, any existing edge between $A_1$ and $V_i$ must be of the form $A_1 \rightarrow V_i$ in  $\g$ and  ${\g}^{*}$.

Next, we show that $\pstar(V_1,Y_1)$ starts with edge $V_1 \leftarrow V_2$ in $\g$.  Since $p$ is chosen as a shortest proper possibly causal path from $A$ to $Y$ that starts with an undirected edge  in $\g$, $p(V_1, Y_1)$ is a proper possibly causal definite status path in $\g$ (\cref{lemma:unshielded-analog-ema-uai}). Then $\pstar(V_1,Y_1)$ is also a path of definite status in $\g^{*}$.  Additionally, since $p(V_1,Y_1)$ is a possibly causal definite status path in $\g$, there cannot be any collider on $\pstar(V_1,Y_1)$.
  
Furthermore, $\pstar$ is a non-causal path, $A_1 - V_1$ is in ${\g}^{*}$, and any edge between $A_1$ and $V_i$, $i \in \{2, \dots , k\}$ is of the form $A_1 \to V_i$, so $\pstar(V_1, Y_1)$ must be a non-causal path from $V_1$ to $Y$. Since  $\pstar(V_1,Y_1)$ is a non-causal definite status path without any colliders, it must start with  an edge into $V_1$, that is $V_1 \leftarrow V_2$ is on $\pstar(V_1,Y_1)$ in $\g^{*}$. 
Then $A_1 - V_1 \leftarrow V_2$ is in ${\g}^{*}$. 

Now, $\pstar(A_1,V_2)$ is of the form $A_1 - V_1 \leftarrow V_2$, so for $\g^{*}$ to be an MPDAG, $\langle A_1, V_2 \rangle $ is in $\g^{*}$ (R1 in \cref{fig:orientationRules}). Then $A_1 \rightarrow V_2 \rightarrow V_1$ and $A_1 - V_1$ are in $\g^{*}$, which by R2 in  Figure \cref{fig:orientationRules} contradicts that $\g^{*}$ is an MPDAG.
\end{proofof}

\begin{proofof}[\cref{thm:dif-effects}]
Statement (iii) directly follows from the construction of the algorithm. Statement (i) follows from the construction of the algorithm and \cref{thm:id-criterion}.

Now we prove statement (ii). The proof follows a similar reasoning as the proof of Theorem 2 of \cite{shpitser2006identification}, proof of Theorem 57 of \cite{perkovic18} and proof of Proposition 3.2.\ of \cite{perkovic20}. 

Suppose for a contradiction that $ |L| \ge 2$ and let ${\g}_{1}$ and $\g_{2}$ be two different MPDAGs in $L$.  Since $\g_{1}$ and $\g_{2}$ are both represented by $\g$, any observational density $f$ consistent with $\g$ is also consistent with $\g_1$ and $\g_2$ due to Markov equivalence.

Let $[\g]$ denote the set of DAGs represented by $\g$. Let $f_1(x_Y|\Do(x_A))$ denote the density of $X_Y$ under the intervention $\Do(X_A = x_A)$ computed from $f(x)$ assuming that the causal DAG belongs to $[\g_{1}]$. Analogously, let $f_2(x_Y|\Do(x_A))$ denote the  density of $X_Y$ under the intervention $\Do(X_A = x_A)$  computed from $f(x)$ assuming that the causal DAG belongs to $[\g_{2}]$.
For the above interventional densities of $X_Y$ to differ, it suffices to show that $\E_1[X_Y| \Do(X_A = \mb{1})] \neq \E_2[X_Y| \Do(X_A = \mb{1})]$, where $\Do(X_A = \mb{1})$ indicates a do intervention that sets the value of every variable indexed by $A$ to 1, and $\E_1$ and $\E_2$ correspond to $f_1$ and $f_2$ respectively. Furthermore, it suffices to show that there is a node $Y_1 \in Y$ such that $\E_1[X_{Y_1}| \Do(X_A = \mb{1})] \neq \E_2[X_{Y_1}| \Do(X_A = \mb{1})]$.

The stages of this proof are as follows. First, we will first establish some graphical differences between $\g_1$ and $\g_2$ that stem from the application of \cref{thm:min-bg} in the \texttt{IDGraphs} algorithm (Algorithm \ref{alg:computegraphs}). These graphical differences will be categorized as cases \ref{case1minimal} and \ref{case2minimal} in Lemma \ref{lemma:auxillaryv1} below.
Then, for each case, we will construct a linear causal model with Gaussian noise that imposes an observational density $f(x)$ consistent with $\g_1$ and $\g_2$ such that $\E_1[X_{Y_1}| \Do(X_A = \mb{1})] \neq \E_2[X_{Y_1}| \Do(X_A = \mb{1})]$, which gives us the desired contradiction.

First, we establish the pertinent graphical differences between $\g_1$ and $\g_2$. 
For this purpose, let $R_1$ and $R_2$ be the list of edge orientations that were added to $\g$ to construct ${\g}_{1}$ and $\g_{2}$ by the \texttt{IDGraphs} algorithm.
That is  $\g_{1}= \MPDAG (\g, R_1)$ and $ \g_{2} = \MPDAG (\g, R_2)$. Without loss of generality, suppose that the edge orientations in $R_1$ and $R_{2}$ are  listed in the order that they were added by the \texttt{IDGraphs} algorithm.

By construction of $R_1$ and $R_2$, there is at least one edge whose orientation differs between $R_1$ and $R_2$. 
Without loss of generality, let $A_1 \rightarrow V_1$, $A_1 \in A$, $V_1 \in V \setminus A$ be the first edge in $R_{1}$ such that $A_1 \leftarrow V_1$ is in $R_2$. 
Also, let $R^{*}$ be the list of edge orientations that come before  $A_1 \rightarrow V_1$ in $R_{1}$ and let $\g^{*}= \MPDAG (\g, R^{*})$. Then by \cref{thm:min-bg}, the total effect of $A$ on $Y$ is not identified given $\g^{*}$.

Among all the shortest proper possibly causal paths from $A$ to $Y$  that start with an undirected edge in $\g^{*}$, choose $\pstar$ as one that starts with $A_1 - V_1$,  $\pstar = \langle A_1 =V_0, V_1, \dots , V_k =Y_1\rangle$, $Y_1 \in Y$.  
Let $p_1$ be the path in a DAG $\g[D]_{1}$ in $[\g_{1}]$ that consists of the same sequence of nodes as $\pstar$ in $\g^{*}$. Analogously, let $p_2$ be  the path in a DAG $\g[D]_{2}$ in $[\g_{2}]$  that consists of the same sequence of nodes as $\pstar$.

By Lemma \ref{lemma:auxillaryv1} we have the following cases:
\begin{enumerate}[label = (\roman*)]
\item \label{case1minimal} if  $\pstar$ is unshielded in $\g^{*}$, then  $p_1$ is of the form $A_1 \rightarrow V_1 \rightarrow \dots \rightarrow Y_1$, and
$p_2$ starts with edge $A_1 \leftarrow V_1$.
\item \label{case2minimal} if $\pstar$  is  a shielded path in $\g^{*}$, then   $A_1 \rightarrow V_i$, $i \in \{1, \dots ,r \}$, $ 2 \le r \le k$, is in $\g^{*}$, $p_2$ is of the form $A_1 \leftarrow V_1 \rightarrow \dots \rightarrow V_r \rightarrow \dots \rightarrow Y_1$, and
\begin{enumerate}
\item \label{case2a} $p_1$ is of the form \\
$A_1 \rightarrow V_1 \rightarrow \dots  \rightarrow V_r \rightarrow \dots \rightarrow Y_1$, or 
\item \label{case2b} $p_1$ is of the form \\
$A_1 \rightarrow V_1 \leftarrow \dots \leftarrow V_l \rightarrow  \dots \rightarrow  Y_1$, $ 2 \le l \le r$.
\end{enumerate}
\end{enumerate}
We will now show how to choose a linear causal model consistent with $\g_1$ and $\g_2$ in each of the above cases that results in $\E_1[X_{Y_1}| \Do(X_A = x_A)] \neq \E_2[X_{Y_1}| \Do(X_A = x_A)]$.

\begin{figure}
\centering
\begin{subfigure}{.45\textwidth}
\centering
\begin{tikzpicture}[>=stealth',shorten >=1pt,auto,node distance=10cm,main node/.style={minimum size=0.4cm,font=\sffamily\Large\bfseries},scale=0.7,transform shape]
\node[main node] (a1) at (0,0) {$A_1$};
\node[main node] (v1) at (1.5,0) {$V_1$};
\node[main node] (n1) at (2.5,0) {};
\node[main node] (n2) at (3.5,0) {};
\node[main node] (y) at (4.5,0) {$Y_1$};

\draw [->,color=red] (a1) edge  (v1);
\draw [->] (v1) edge  (n1);
\draw [dotted] (n1) edge (n2);
\draw [->] (n2) edge (y);
\end{tikzpicture}
\caption{}
\end{subfigure}
\vrule
\begin{subfigure}{.45\textwidth}
\centering
\begin{tikzpicture}[>=stealth',shorten >=1pt,auto,node distance=10cm,main node/.style={minimum size=0.4cm,font=\sffamily\Large\bfseries},scale=0.7,transform shape]
\node[main node] (a1) at (0,0) {$A_1$};
\node[main node] (v1) at (1.5,0) {$V_1$};
\node[main node] (n1) at (2.5,0) {};
\node[main node] (n2) at (3.5,0) {};
\node[main node] (y) at (4.5,0) {$Y_1$};

\draw [->,color=red] (v1) edge  (a1);
\draw [->] (v1) edge  (n1);
\draw [dotted] (n1) edge (n2);
\draw [->] (n2) edge (y);
\end{tikzpicture}
\caption{}
\end{subfigure}
\caption{DAGs (a) $\g[D]_{11}$ and (b) $\g[D]_{21}$ corresponding to \ref{case1minimal} in the Proof of \cref{thm:dif-effects}.}
\label{fig:case1}
\end{figure}

\ref{case1minimal} Consider a multivariate Gaussian density over $X$ with mean zero, constructed using a linear causal model with Gaussian noise consistent with $\g[D]_1$ and thus, also $\g_1$ (due to Markov equivalence). 
We define the linear causal model  in such a way that all edge coefficients except for the ones on $p_1$ are $0$, and all edge coefficients  on $p_1$ are in $(0,1)$ and small enough so that we can choose the error variances in such a way that $\Var(X_i)= 1$ for every $i \in V$.

The density $f(x)$ generated in this way is consistent with $\g[D]_{1}$ and thus also consistent with $\g_1$ and $\g_{2}$ \citep{lauritzen1990independence}. Moreover, $f(x)$ is  consistent with  DAG $\g[D]_{11}$ that is obtained from $\g[D]_{1}$ by removing all edges except for the ones on $p_1$; see Figure \ref{fig:case1}(a). Additionally, $\g[D]_{11}$ is Markov equivalent to DAG $\g[D]_{21}$, which is obtained from $\g[D]_{2}$ by removing all edges except for those on $p_2$; see Figure \ref{fig:case1}(b). Hence, $f(x)$ is also consistent with $\g[D]_{21}$. 

Let  $f_{1}(x_{Y_1} |\Do(x_{A}))$ be an interventional density of $X_{Y_1}$ under the intervention $\Do(X_A = x_A)$ that is consistent with $\g[D]_{11}$ (and  $\g[D]_1$).
By Rules 3 and 2 of the do-calculus (\cref{lem:do-rules}), we have
\begin{align*}
f_1(x_Y | \Do(x_A)) =  f_1 (x_Y| \Do(x_{A_1})) = f(x_Y |x_{A_1}).
\end{align*} 
So $\E_1[X_{Y_1}|\Do(X_A = \mb{1})] =\int x_{Y_1} f(x_Y| X_{A_1} =\mb{1}) \dd x_{Y_1} =  \Cov(X_{Y_1}, X_{A_1}) = a$ by \cref{lem:mardia-condexp}. Additionally, by \cref{lem:wright},  $a$ is equal to the product of all edge coefficients along $p_1$ and so $a \in (0, 1)$. 

Similarly, let $f_{2}(x_{Y_1} |\Do(x_{A}))$ be an interventional density of $X_{Y_1}$ consistent with $\g[D]_{21}$ (and $\g[D]_{2}$).
Then  $f_{2}(x_{Y_1} |\Do(x_{A})) =f(x_{Y_1})$ by Rule 3 of \cref{lem:do-rules}. Hence, $\E_2[X_{Y_1}|\Do(X_A = \mb{1} )] = \E[X_{Y_1}] =  0$. Since $a \neq 0$, this completes the proof for case \ref{case1minimal}.

\begin{figure}
\centering
\begin{subfigure}{.48\textwidth}
\centering
\begin{tikzpicture}[>=stealth',shorten >=1pt,auto,node distance=10cm,main node/.style={minimum size=0.4cm,font=\sffamily\Large\bfseries},scale=0.7,transform shape]
\node[main node] (a1) at (0,0) {$A_1$};
\node[main node] (v1) at (1.2,0) {$V_1$};
\node[main node] (v2) at (2.4,0) {$V_2$};
\node[main node] (n1) at (3.6,0) {\normalsize $\cdots$};
\node[main node] (vr) at (4.8,0) {$V_r$};
\node[main node] (n3) at (6,0) {\normalsize $\cdots$};
\node[main node] (y) at (7.2,0) {$Y_1$};

\node[main node] (DOT1) at (3.5,0.6) {\normalsize $\cdots$};

\draw [->,color=red] (v1) edge  (a1);
\draw [->] (v1) edge  (v2);
\draw [->] (v2) edge  (n1);
\draw [->] (n1) edge (vr);
\draw [->] (vr) edge  (n3);
\draw [->] (n3) edge (y);
\draw[->] (a1) ..  controls (.6, 1.5) and (1.5, 1.5) ..  (v2);
\draw[->] (a1) ..  controls (0.4, 2) and (4, 2) ..  (vr);
\end{tikzpicture}
\caption{}
\end{subfigure}
\begin{subfigure}{.48\textwidth}
\centering
\begin{tikzpicture}[>=stealth',shorten >=1pt,auto,node distance=10cm,main node/.style={minimum size=0.4cm,font=\sffamily\Large\bfseries},scale=0.7,transform shape]
\node[main node] (a1) at (0,0) {$A_1$};
\node[main node] (v1) at (1.2,0) {$V_1$};
\node[main node] (v2) at (2.4,0) {$V_2$};
\node[main node] (n1) at (3.6,0) {\normalsize $\cdots$};
\node[main node] (vr) at (4.8,0) {$V_r$};
\node[main node] (n3) at (6,0) {\normalsize $\cdots$};
\node[main node] (y) at (7.2,0) {$Y_1$};

\node[main node] (DOT1) at (3.5,0.6) {\normalsize $\cdots$};

\draw [->,color=red] (a1) edge  (v1);
\draw [->] (v1) edge  (v2);
\draw [->] (v2) edge  (n1);
\draw [->] (n1) edge (vr);
\draw [->] (vr) edge  (n3);
\draw [->] (n3) edge (y);
\draw[->] (a1) ..  controls (.6, 1.5) and (1.5, 1.5) ..  (v2);
\draw[->] (a1) ..  controls (0.4, 2) and (4, 2) ..  (vr);
\end{tikzpicture}
\caption{}
\end{subfigure}
\begin{subfigure}{\columnwidth}
\centering
\begin{tikzpicture}[>=stealth',shorten >=1pt,auto,node distance=10cm,main node/.style={minimum size=0.4cm,font=\sffamily\Large\bfseries},scale=0.65,transform shape]
\node[main node] (a1) at (0,0) {$A_1$};
\node[main node] (v1) at (1.5,0) {$V_1$};
\node[main node] (v2) at (3,0) {$V_2$};
\node[main node] (n1) at (4,0) {};
\node[main node] (n2) at (5,0) {};
\node[main node] (vl) at (6,0) {$V_l$};
\node[main node] (n3) at (7,0) {};
\node[main node] (n4) at (8,0) {};
\node[main node] (vr) at (9,0) {$V_r$};
\node[main node] (n5) at (10,0) {};
\node[main node] (n6) at (11,0) {};
\node[main node] (y) at (12,0) {$Y_1$};

\node[main node] (DOT1) at (3.8,0.7) {\normalsize $\cdots$};
\node[main node] (DOT2) at (6.5,0.7) {\normalsize $\cdots$};

\draw [->,color=red] (a1) edge  (v1);
\draw [->,color=red] (v2) edge  (v1);
\draw [->,color=red] (n1) edge  (v2);
\draw [dotted] (n1) edge (n2);
\draw [->,color=red] (vl) edge (n2);
\draw [->] (vl) edge  (n3);
\draw [dotted] (n3) edge (n4);
\draw [->] (n4) edge (vr);
\draw [->] (vr) edge  (n5);
\draw [dotted] (n5) edge (n6);
\draw [->] (n6) edge (y);
\draw[->] (a1) ..  controls (1.2, 1.5) and (2, 1.5) ..  (v2);
\draw[->] (a1) ..  controls (1.0, 1.9) and (4, 1.9) ..  (vl);
\draw[->] (a1) ..  controls (0.8, 2.2) and (6, 2.2) ..  (vr);
\end{tikzpicture}
\caption{}
\end{subfigure}
\caption{DAGs (a) $\g[D]_{22}$, (b) $\g[D]_{11}$, and (c) $\g[D]_{12}$ corresponding to \ref{case2a} and \ref{case2b} in the Proof of \ref{thm:dif-effects}.}
\label{fig:case2a}
\end{figure}

\ref{case2minimal} Consider a multivariate Gaussian density over $X$ with mean zero, constructed using a linear causal model  with Gaussian noise consistent with $\g[D]_2$. 
We define the causal model  in  a way such that all edge coefficients except for the ones on $p_2$ and $A_1 \to V_i$, $i \in \{2, \dots , r \}$ are $0$, and all edge coefficients  on $p_2$ and $A_1 \to V_i$ are in $(0,1)$ and small enough so in such a way that $\Var(X_i)= 1$ for all $i \in V$. 

The density $f(x)$ generated in this way is consistent with  $\g[D]_2$ and $\g_{2}$ \citep{lauritzen1990independence}. Moreover, $f(x)$ is  consistent with  DAG $\g[D]_{22}$ that is obtained from $\g[D]_{2}$ by removing all edges except for the ones on $p_2$  and $A_1 \to V_i$, $i \in \{2, \dots , r \}$, $2 \le r \le k$; see \cref{fig:case2a}(a).
Let  $f_{2}(x_{Y_1} |\Do(x_{A}))$ be an interventional density of $X_{Y_1}$ under the intervention $\Do(X_A = x_A)$ that is consistent with $\g[D]_{22}$ (and also $\g[D]_2$). 

We now have
\begin{align}
f_2(x_{Y_1} | \Do(x_A)) &= f_2(x_{Y_1}| \Do(x_{A_1}))\nonumber \\
& = \int f(x_{Y_1}| \Do(x_{A_1}), x_{V_1}) f(x_{V_1} | \Do(x_{A_1})) \dd x_{V_1} \nonumber\\
&= \int f(x_{Y_1} |x_{A_1}, x_{V_1}) f(x_{V_1}) \dd x_{V_1}. \label{eq:case21} 
\end{align}
The first line follows using Rule 3 of the do-calculus, and the third line follows from an application of Rule 2 and Rule 3; see \cref{lem:do-rules}.

We now compute $\E_{2}[X_{Y_1}| \Do(X_A = \mb{1})] $. For simplicity, we will use shorthands $ \Cov(X_{Y_1}, X_{A_1}) = a$,  $\Cov(X_{Y_1}, X_{V_1}) = b$ and $\Cov(X_{A_1},X_{V_1}) =c$.
Now, using \cref{lem:mardia-condexp} and \cref{eq:case21}, we have 
\begin{align*}
\E_2[X_{Y_1}| \Do(X_A = 1)] &= \int \E[X_{Y_1}| X_{A_1} = 1, X_{V_1} =x_{V_1}] f(x_{V_1}) \dd x_{V_1}\\
&= \int \begin{bmatrix}
       a & b
     \end{bmatrix}
     \begin{bmatrix}
       1 & c\\
     c & 1
     \end{bmatrix}^{-1}
     \begin{bmatrix}
       1  \\
       x_{V_1}
     \end{bmatrix} f(x_{V_1}) \dd x_{V1}\\
&=\int  \frac{1}{1- c^2}\begin{bmatrix}
       a&     b
     \end{bmatrix}
     \begin{bmatrix}
       1 & -c \\
       -c & 1
     \end{bmatrix}
     \begin{bmatrix}
       1  \\
       x_{V_1}
     \end{bmatrix} f(x_{V_1}) \dd x_{V1}\\
&=\frac{a - bc}{1-c^2} + \frac{- ac + b}{1-c^2} \E[X_{V_1}]= \frac{a - bc}{1-c^2} .
\end{align*}

Now, consider the cases \ref{case2a} and \ref{case2b}. Note that $f(x)$ is also consistent with 
 $\g[D]_1$, $\g_1$ and a $\DAG$ that is  obtained from $\g[D]_{1}$ by removing all edges except for the ones on $p_1$  and $A_1 \to V_i$, $i \in \{2 ,\dots , r\}$ \citep{lauritzen1990independence}. Depending on case  \ref{case2a} or \ref{case2b}, this will be either DAG $\g[D]_{11}$ in Figure \ref{fig:case2a}(b) or DAG $\g[D]_{12}$ in  Figure \ref{fig:case2a}(c). 

Let $f_{11}(x_{Y_1} |\Do(x_{A})) $ and $f_{12}(x_{Y_1} |\Do(x_{A}))$ be the interventional densities of $X_{Y_1}$ that are consistent with $\g[D]_{11}$ and $\g[D]_{12}$, respectively. 
Note that  $f_{11}(x_{Y_1} |\Do(x_{A})) = f_{12}(x_{Y_1} |\Do(x_{A}))$ since
\begin{align}
f_{11}(x_{Y_1}|\Do(x_{A}) ) &= f_{11}(x_{Y_1} | \Do(x_{A_1})) = f(x_{Y_1} | x_{A_1}) \nonumber\\
&= f(x_{Y_1} | x_{A})= f_{12} (x_{Y_1} |\Do(x_A)). \label{eq:case2ab}
\end{align}
The first two equalities above follow from Rule 3 and Rule 2 of the do-calculus, while the third and forth follow from Rule 1 and Rule 2; see again \cref{lem:do-rules}. Hence, $f_1(x_{Y_1} |\Do(x_A)) = f_{11}(x_{Y_1} |\Do(x_{A})) = f_{12}(x_{Y_1} |\Do(x_{A}))$. 

Using Equation \eqref{eq:case2ab} and \cref{lem:mardia-condexp}, we have $\E_1[X_{Y_1}| \Do(X_A = \mb{1})] = \E[X_{Y_1} | X_{A_1} = \mb{1}] = \Cov(X_{Y_1}, X_{A_1}) = a$. To show that $\E_1[X_{Y_1} | \Do(X_A = \mb{1}) ] \neq \E_2[X_{Y_1} | \Do(X_A =1)]$, we need only to show that $a \neq (a - bc)/(1-c^2)$. 

We will show that $b > ac$ and $c >0$, which leads to $a - bc < a - ac^2$, that is  $(a - bc)/(1-c^2) < a$. 
To show $b > ac$ and $c>0$, we need to discuss $a, b$ and $c $  in terms of the original linear causal model. 

By \cref{lem:wright}, we have that $c = \Cov(X_{A_1}, X_{V_1})$ is equal the edge coefficient assigned to $A_1 \leftarrow V_1$ in $\g[D]_{21}$, and hence $c \in (0,1)$.  Let $a_1$ be the product of edge coefficients on $p_2(V_1, Y_1)$ and let $a_i$ be the product of edge coefficients along  $\langle A_1 , V_i \rangle  \oplus p_2(V_i, Y_1)$, $ i \in \{2, \dots , r \}$.  Then $a_i \in (0,1)$ for all $i \in \{1, \dots , r\}.$
By \cref{lem:wright}, we now have
\begin{align*}
&a = \Cov(X_{Y_1}, X_{A_1}) = c\cdot a_1 + a_2 + \cdots + a_r,\\
&b = \Cov(X_{Y_1}, X_{V_1}) = a_1 + c \cdot (a_2 + \cdots + a_r),
\end{align*} 
which yields $b - ac = a_1 (1 - c^2) > 0$, completing the proof.
\end{proofof}

\begin{lemma}\label{lemma:auxillaryv1}
Suppose that the total effect of $A$ on $Y$ is not identified given $\MPDAG$ $\g$.  Let $ p = \langle A_1= V_0,  V_1, \dots, V_k = Y_1 \rangle$, $k \ge 1$, $A_1 \in A$, $Y_1 \in Y$, be a shortest proper possibly causal path from $A$ to $Y$ in $\g$.
Let ${\g}_{1} = \MPDAG(\g, \{A_1 \rightarrow V_1\})$ and ${\g}_{2} =  \MPDAG({\g,  \{A_1 \leftarrow V_1\}})$. Let $p_1$ and $p_2$ be the paths in $\g_{1}$ and $\g_{2}$ respectively, that consist of the same sequence of nodes as $p$ in $\g$.

\begin{enumerate}[label = (\roman*)]
\item\label{cond-unshielded} If $p$ is an unshielded path in $\g$, then
\begin{itemize}
\item $p_1$ is of the form $A_1 \rightarrow V_1 \rightarrow \dots \rightarrow Y_1$, and
\item $p_2$ is of the form $A_1 \leftarrow V_1 \dots Y_1$.
\end{itemize}
\item\label{cond-shielded} If $p$ is a shielded path in $\g$, then 
\begin{itemize}
\item $A_1 \rightarrow V_i$ is in $\g$ for all $i \in \{2, \dots , r\}, r \le k$, $k >1$,
\item $p_2$ is of the form $A_1 \leftarrow V_1 \rightarrow \dots \rightarrow Y$, 
\item Let $\g[D]_{1}$ be a DAG in $[\g_{1}]$ and let $p_{11}$ be the path in $\g[D]_{1}$ corresponding to $p_1$ in $\g_1$ and to $p$ in $\g$, then 
\begin{enumerate}
\item $p_{11}$ is of the form $A_1 \rightarrow V_1 \rightarrow \dots \rightarrow Y_1$ in $\g[D]_{1}$, or
\item $p_{11}$ is of the form $A_1 \rightarrow V_1 \leftarrow \dots \leftarrow V_l \rightarrow \dots  \rightarrow Y_1$, $1< l \le r$ in $\g[D]_{1}$.
\end{enumerate}
\end{itemize}
\end{enumerate}
\end{lemma}

\begin{proofof}[Lemma \ref{lemma:auxillaryv1}]
Path $p$ is chosen as a shortest proper possibly causal path from $A$ to $Y$ that starts with an undirected edge in $\g$. Hence, $p(V_1, Y_1)$ must be an unshielded possibly causal path from $V_1$ to $Y_1$, otherwise we can choose a shorter path than $p$ in $\g$. 
This implies that no node $V_i$, $i \in \{2, \dots k-1\}$ can be a collider on either $p_1$ or $p_2$.

\ref{cond-unshielded} Suppose first that $p$ itself is unshielded. That is,  no edge  $\langle V_{i}, V_{i+2}\rangle$, $i \in \{ 0, k-2\}$ is in $\g$.
Of course, since $\g_2$ contains edge $A_1 \leftarrow V_1$, $p_2$ is of the form $A_1 \leftarrow V_1 \dots Y_1$. Hence, we only need to show $p_1$ is a causal path in $\g_2$.

Since $p$ is unshielded, $p_1$ is also an unshielded path.
Since $A_1 \rightarrow V_1$ is in $\g_{1}$, as a consequence of iterative application of rule R1 of \cite{meek1995causal} (\cref{fig:orientationRules}), $p_1$ is a causal path in $\g_{1}$. 

\ref{cond-shielded}
Next, we suppose that $p$ is shielded. We first show that $A_1 \rightarrow V_i$, for all $i \in \{2, \dots , r\}$, $r \le k$ is in $\g$.

  As discussed at the beginning of this proof,  $p(V_1, Y_1)$ is unshielded. Therefore, since $p$ is shielded and $p(V_1, Y_1)$ is unshielded,  edge $\langle A_1  ,V_2 \rangle $ is in $\g$.  Furthermore, since $p$ is chosen as a shortest proper possibly causal path from $A$ to $Y$ that starts with an undirected edge in $\g$, $\langle A_1, V_2 \rangle$ must be of the form $A_1 \rightarrow V_2$. 

If path $\langle A_1, V_2 \rangle \oplus p(V_2, Y_1)$ is shielded, then by the same reasoning as above, $A_1 \rightarrow V_3$ is in $\g$. We can continue with the same reasoning, until we reach $V_r$, $r \in \{2, \dots , k \}$, so that  $A_1 \rightarrow V_i$ is in $\g$ for $i \in \{2, \dots , r\}$ and  $\langle A_1, V_r \rangle \oplus p(V_r, Y_1)$ is an unshielded possibly causal path. 

We note that if $r<k$,  $p(V_r, Y_1)$ is of the form $V_r \rightarrow \dots \rightarrow Y_1$. This is due to the fact that  $A_1 \rightarrow V_r$ is in $\g$ and that $\langle A_1, V_r \rangle \oplus p(V_r, Y_1)$ is an unshielded possibly causal path in $\g$.

Next, we show that $p_2$ is of the form $A_1 \leftarrow V_1 \rightarrow \dots \rightarrow Y_1$.
 Since $A_1 \leftarrow V_1$ and $A_1 \rightarrow V_2$ are in $\g_{2}$ and since $\g_{2}$ is acyclic, by rule R2 of \cite{meek1995causal}, the edge $\langle V_1, V_2 \rangle$ is of the form $V_1 \rightarrow V_2$ in $\g_2$.
Then since $p_2(V_1, Y_1)$ is an unshielded possibly causal path that starts with $V_1 \rightarrow V_2$, by iterative applications of rule R1 of \cite{meek1995causal}, $p_2(V_1, Y_1)$ must be a causal path in $\g_{2}$. 

Suppose that $\g[D]_1$ is a DAG in $[\g_1]$. Based on the above, we know that $A_1 \rightarrow V_1$ and if $r <k$,  $V_r \rightarrow \dots Y_1$ are in $\g_1$ and therefore, in $\g[D]_1$ as well.
The subpath $p(V_1, V_r)$ is a possibly causal unshielded path in $\g$ and hence, no node among $V_2, \dots , V_{r-1}$ is a collider on $p$, $p_1$, or $p_{11}$. 
Therefore, either $p_{11}(V_1, V_r)$ is a causal path in $\g[D]_{11}$, in which case $p_{11}$ is of the form $A_1 \rightarrow V_1 \rightarrow \dots \rightarrow Y_1$, or there is a node $V_l$, $1 < l \le r$ on $p_{11}(V_{11}, V_r)$, such that, $p_{11}(V_1, V_r)$ is of the form $V_1 \leftarrow \dots \leftarrow  V_l \rightarrow \dots \rightarrow V_r$.
\end{proofof}

\bibliographystyle{apalike}

\end{document}